\documentclass[a4paper,10pt]{amsart}
\usepackage{amsmath,amsthm,latexsym,amscd,amssymb}
\usepackage[all]{xy}
\usepackage{colortbl}
\usepackage{hyperref}
\numberwithin{equation}{section}

\newtheorem{prop}{Proposition}[section]
\newtheorem{lem}[prop]{Lemma}
\newtheorem{dfn}[prop]{Definition}
\newtheorem{theorem}[prop]{Theorem}

\newtheorem{rem}[prop]{Remark}
{\par \vspace{0cm}\ \\ \noindent \%\color{blue}\% \begin{tiny}}{\par\end{tiny}\%\color{black}\%\vspace{0cm}}
{\par \vspace{0cm}\ \\ \noindent \%\color{red}\% \begin{tiny}}{\par\end{tiny}\%\color{black}\%\vspace{0cm}}

\def\id\textrm{Id}
\def\dd{\mathcal{D}(d)}

\def\kerd{\operatorname{ker}(d)}

\def \dl{\partial L_x^0}

\def\t{\tau}

\def\D{\Delta}

\def\R{\mathbin{\mathbb R}}

\def\C{\mathbb{C}}

\newcommand{\ra}{\operatorname{range}}

\newcommand{\clis}{H^{2k}_{(2),dR}(L_x^0)}
\newcommand{\cali}{L^2(\Omega^{2k}(\partial L_x^0))}

\newcommand{\ty}{\infty}
\definecolor{light}{gray}{.95}

\newcommand\Di{D\kern-7pt/}

\title{The Atiyah Patodi Singer signature formula for measured foliations}

\author{Paolo Antonini}
\address{Paolo Antonini\\
Universit\"at Regensburg\\   Universit\"at{\ss}tra{\ss}e 31  93053 Regensburg
\\
Germany}
\email{paolo.anton@gmail.com, paolo.antonini@mathematik.uni-regensburg.de}

\thanks{}

\begin{document}
\begin{abstract}
\noindent Let $(X_0,\mathcal{F}_0) $ be a compact manifold with boundary endowed with a  foliation 
$\mathcal{F}_0$ which is assumed to be measured and transverse  to the boundary. 
We denote by $\Lambda$  a holonomy invariant transverse measure on  $(X_0,\mathcal{F}_0) $ and by 
$\mathcal{R}_0$ the equivalence relation of the foliation.
Let $(X,\mathcal{F})$ be the corresponding manifold with 
 cylindrical ends  and extended foliation with equivalence relation $\mathcal{R}$. In the paper \cite{io} we proved a  formula for the $L^2$-$\Lambda$ index of a 
 longitudinal 
 Dirac-type operator $D^{\mathcal{F}}$ on $X$ in the spirit of Alain Connes' non commutative geometry  \cite{Cos}.
 Here we specialize ourselves to the signature operator. 
\noindent We define three  types of signature for the pair (foliation, fboundary foliation): the analytic signature, denoted
 $\sigma_{\Lambda,\operatorname{an}}(X_0,\partial X_0)$ is the leafwise $L^2$-$\Lambda$-index of the signature operator on the extended manifold;
 the Hodge signature 
 
 \noindent $\sigma_{\Lambda,\operatorname{Hodge}}(X_0,\partial X_0)$, defined using the natural representation of $\mathcal{R}$ on the field of square integrable harmonic forms on the leaves and
the de Rham signature, $\sigma_{\Lambda,\operatorname{dR}}(X_0,\partial X_0)$, defined using the natural representation  
of $\mathcal{R}_0$ on the field of the $L^2$  relative de Rham spaces of the leaves.
We prove that 
these three signatures coincide 
$$\sigma_{\Lambda, \operatorname{an}}(X_0,\partial X_0)=
\sigma_{\Lambda,\operatorname{Hodge}}(X_0,\partial X_0)=\sigma_{\Lambda,\operatorname{dR}}(X_0,\partial X_0).$$ 
As a consequence of these equalities and of the index formula we finally obtain the main result of this work,
the Atiyah-Patodi-Singer
signature formula for measured foliations:
$$\sigma_{\Lambda,\operatorname{dR}}(X_0,\partial X_0)=\langle L(T\mathcal{F}_0),C_{\Lambda}\rangle +1/2[\eta_{\Lambda}(D^{\mathcal{F}_{\partial}})].$$ We give also, in the appendix an account of Non commutative integration theory.
\end{abstract}
\maketitle

\tableofcontents

\section{Introduction}
Let $X_0$ be a $4k$--dimensional oriented manifold without boundary. One can give four different definitions of the signature.
\begin{itemize}
\item The \underline{topological signature} $\sigma(X_0)$ is defined as the signature of the intersection form in the middle degree cohomology;
$(x,y):=\langle x \cup y,[X_0]\rangle, \,x,y \in \operatorname{H}^{2k}(X_0,\R).$
\item The \underline{de Rham signature} $\sigma_{\operatorname{dR}}(X_0)$ is the signature of the Poincar\'e intersection form in the middle de Rham cohomology;
$([\omega],[\phi]):=\int_{X_0}\omega\wedge \phi;\,\,\omega,\phi \in \operatorname{H}^{2k}_{dR}(X_0).$
\item The \underline{Hodge signature}, $\sigma_{\operatorname{Hodge}}(X_0)$ is the signature of the Poincar\'e intersection form defined in the space of $2k$ Harmonic forms with respect to some Riemannian metric;
$(\omega,\phi):=\int_{X_0}\omega\wedge \phi;\,\,\omega,\phi \in \mathcal{H}^{2k}(X_0).$
\item The \underline{analytical signature} is the index of the chiral signature operator\footnote{this is the differential operator $d+d^*$ acting on the complex of differential forms, odd w.r.t. the natural chiral grading $\tau:=(-1)^k\ast (-1)^{\frac{|\cdot|(|\cdot|-1)}{2}}$, $D^{\operatorname{sign}}=\left(\begin{array}{cc}0 & D^{\operatorname{sign},-} \\D^{\operatorname{sign},+} & 0\end{array}\right)$}
$$\sigma_{\operatorname{an}}(X_0):=\operatorname{ind}(D^{\operatorname{sign},+}).$$
\end{itemize} One can prove that all these numbers coincide,
\begin{equation}\label{allequal}
\sigma(X_0)=\sigma_{\operatorname{dR}}(X_0)=\sigma_{\operatorname{Hodge}}(X_0)=\sigma_{\operatorname{an}}(X_0).
\end{equation} 
 The Hirzebruch formula $$\sigma(X_0)=\int_{X_0}L(X_0)$$ can be proven using cobordism arguments as in the original work of Hirzebruch or can be seen as a consequence of the Atiyah--Singer index formula together with Hodge theory \cite{BeGeVe}.

\noindent If $\widetilde{X_0}\longrightarrow X_0$ is a Galois covering with deck group $\Gamma$ Atiyah \cite{At-e} used the Von Neumann algebra of the group $\Gamma$ to normalize the signature on the $L^2$--middle degree harmonic forms on the total space. This signature $\sigma_{\Gamma}(\widetilde{X_0})$ again  enters in a Hirzebruch type formula
$$\sigma_{\Gamma}(\widetilde{X_0})=\int_{X_0}L(X_0)$$
hence turns out to be equal to the signature of the base. This is the celebrated Atiyah $L^2$--signature theorem.

\noindent The Atiyah $L^2$--signature theorem was extended by Alain Connes \cite{Cos} to the situation in which $X_0$ is foliated by an even dimensional foliation. This is the realm of the \emph{non--commutative geometry}. 

\noindent
What can one say if $X_0$ has non empty boundary ?

\noindent So let now $X_0$ be an oriented compact manifold with boundary and suppose the metric is product type near the boundary. Attach an infinite cylinder across the boundary to form the manifold with cylindrical ends $X=X_0\bigcup_{\partial X_0}\Big{[} \partial X_0\times [0,\infty)\Big{]}$.
In the seminal paper by Atiyah Patodi and Singer \cite{AtPaSi1} is showen that the  Fredholm index of the generalized boundary value problem
with the pseudodifferential A.P.S. boundary condition
 on $X_0$ for the signature operator is connected to the $L^2$--index of the extended operator on $X$. Indeed this Fredholm index is the $L^2$--index on $X$ plus a defect related to the space of extended solutions on the cylinder. More precisely the operator on the cylinder acting on the natural space of $L^2$--sections is no more Fredholm (in the general case in which the boundary operator is not invertible) but its kernel and the kernel of its formal adjoint are finite dimensional and the difference of the dimensions is given  by the formula\footnote{opposite orientation w.r.t. A.P.S.} $$\operatorname{ind}_{L^2}(D^+)=\int_{X_0} \hat{A}(X,\nabla)\operatorname{Ch}(E)+\dfrac{\eta(0)}{2}+\dfrac{h_{\infty}(D^-)-h_{\infty}(D^+)}{2};$$
  where $h_\infty(D^{\pm})$ are the dimensions of the limiting values of the extended $L^2$--solutions and  $\eta(0)$ is the eta invariant of the boundary operator. 
Then, in the case of the signature operator the authors investigate the relationship between the A.P.S. index of the operator on $X_0$, the (topological) signature of the pair $(X_0,\partial X_0)$, the $L^2$--index on $X$ and the space of square integrable harmonic forms on $X$.
The result of A.P.S. says that the signature $\sigma(X_0)$ is exactly the $L^2$--index on the cylinder i.e. the difference of the dimensions $h^\pm$ of positive/negative square integrable harmonic forms\footnote{indeed the intersection form passes to be non--degenerate to the image of the relative cohomology into the absolute one. This vector space is naturally isomorphic to the space of $L^2$--harmonic forms on $X$.} 
on $X$,
$\sigma(X_0)=h^+-h^-=\operatorname{ind}_{L^2}(D^{\operatorname{sign},+})$ while $h_{\infty}(D^{\operatorname{sign},-})=h_{\infty}(D^{\operatorname{sign},+})$ by specific simmetries of the signature operator.
In particular the A.P.S. signature formula becomes 
$$\sigma(X_0)=\int_{X_0}L(X_0,\nabla)+\eta\big{(}D^{\operatorname{sign}}_{|\partial X_0}\big{)}.$$

\noindent In the case of $\Gamma$--Galois coverings of a manifold with boundary
 with a cylinder attached, $\widetilde{X}\longrightarrow X$ this program is partially carried out by Vaillant \cite{Vai} in his Master thesis. Vaillant estabilishes a Von Neumann index formula in the sense of Atiyah \cite{At-e} for a Dirac type operator
 and relates this index with the $\Gamma$--dimensions of the harmonic forms on $\widetilde{X}$. The remaining part of the story i.e. the relation with the topologically defined $L^2$--signature is developed by L{\"u}ck and Schick \cite{lus}. 
Call the index of Vaillant the analytical $L^2$--signature of the covering $\widetilde{X_0}\rightarrow X_0$ over the compact piece $X_0$, in symbols $\sigma_{\operatorname{an},(2)}(\widetilde{X_0})$ while $\sigma_{\textrm{Hodge}}(\widetilde{X_0})$ is the $L^2$ signature defined using harmonic forms on $\widetilde{X}$. Vaillant proves 
 $$\sigma_{\textrm{an},(2)}(\widetilde{X_0})=\int_{X_0}L(X_0,\nabla)+\eta_{\Gamma}\big{(}D^{\operatorname{sign}}_{|\partial \tilde{X}}\big{)}=\sigma_{\operatorname{Hodge}}(\widetilde{X_0}).$$
Luck and Schick define other different types of $L^2$--signatures; de Rham $\sigma_{\operatorname{dR},(2)}(\widetilde{X_0})$ and simplicial $\sigma_{\operatorname{top},(2)}(\widetilde{X_0})$ and prove that they are all the same and coincide with the signatures of Vaillant
$\sigma_{\operatorname{Hodge}}(\widetilde{X_0})=\sigma_{\operatorname{dR},(2)}(\widetilde{X_0})=\sigma_{\operatorname{top},(2)}(\widetilde{X_0}).$
None of these steps are easy adaptations of the closed case since in the classical proof a fundamental role is played by the existence of a gap around the zero in the spectrum of the boundary operator. This situation fails to be true in non compact (or cocompact) ambients.

\bigskip
\noindent I this paper we carry out this program for a foliated manifold with cylindrical ends  endowed with a holonomy invariant measure $\Lambda$ \cite{Cos}. The framework is that defined by Connes in his seminal paper about non commutative integration theory \cite{Cos}. Making use in a crucial way of the various semifinite Von Neumann algebras associated to square integrable representations of the Borel groupoid of the foliation (its equivalence relation $\mathcal{R}$) the author extended  the index formula of Vaillant \cite{io}.
  \begin{theorem}
The Dirac operator has finite $L^2$--$\Lambda$--index and the following formula holds
\begin{align}\label{2111}\operatorname{ind}_{L^2,\Lambda}(D^+)=\langle\widehat{A}(X)\operatorname{Ch}(E/S),C_{\Lambda}\rangle +1/2[\eta_{\Lambda}(D^{\mathcal{F}_{\partial}})-h^+_{\Lambda}+h^{-}_{\Lambda}].\end{align} 
\end{theorem}
\noindent The dimensions of the spaces of extended solutions, $h^\pm_{\Lambda}$  are suitably defined using the fields of extended solutions along the leaves. The foliation eta invariant is defined by Ramachandran \cite{Rama} and the usual integral in the A.P.S. formula is changed into the  distributional pairing of a tangential distributional form with the Ruelle--Sullivan current \cite{MoSc}. Some work is needed to show that for the signature operator $h^+_{\Lambda}=h^-_{\Lambda}$. Inspired by the definitions of L{\"u}ck and Schick \cite{lus} we pass study three different representations of $\mathcal{R}_0$ (the equivalence relation of the foliation on the compact piece $X_0$) in order to define the \underline{Analytical Signature}, $\sigma_{\Lambda,\operatorname{an}}(X_0,\partial X_0)$
(i.e. the measured index of the signature operator on the cylinder), the \underline{de Rham signature} 
$\sigma_{\Lambda,\operatorname{dR}}(X_0,\partial X_0)$
(i.e the one induced by the representation which is valued in the relative de Rham spaces of the leaves) and the \underline{Hodge signature}, $\sigma_{\Lambda,\operatorname{Hodge}}(X_0,\partial X_0)$ (defined in terms of the representation of $\mathcal{R}_0$ in the harmonic forms on the leaves of the foliation on $X$).

\noindent Combining a generalization of the notion of the $L^2$ long exact sequence of the pair $(\mathcal{F}_0,\partial \mathcal{F})$, in the sense of sequences of Random Hilbert complexes (the analog of the homology $L^2$ long sequence of Hilbert $\Gamma$--modules in Cheeger and Gromov \cite{ChGr}) together with the analysis of boundary value problems of \cite{sch2}, we show that the methods in \cite{lus} can be generalized to give the following

\begin{theorem}
The above three notions of the $\Lambda$--signature for the foliation on $X_0$ coincide,
$$\sigma_{\Lambda,\operatorname{dR}}(X,\partial X_0)=\sigma_{\Lambda,\operatorname{Hodge}}(X,\partial X_0)=\sigma_{\Lambda,\operatorname{an}}(X,\partial X_0)$$
 and the following A.P.S. signature formula holds true
$$\sigma_{\Lambda, \operatorname{an}}(X_0,\partial X_0)=\langle L(X),C_{\Lambda}\rangle +1/2[\eta_{\Lambda}(D^{\mathcal{F}_{\partial}})].$$
\end{theorem}
\noindent The author wishes to thank Paolo Piazza for having suggested him the problem and for a number of interesting discussions, Moulay T. Benameur, Georges Skandalis, James Heitch, Eric Leichtnam, Stephane Vassout and Yuri Kordyukov for discussions and comments. 
\section{Geometric Setting}\label{geom}
\noindent A $p$--dimensional foliation $\mathcal{F}_0$ on a manifold with boundary $ X_0$ is transverse to the boundary if it is given by a foliated atlas $\{U_{\alpha}\}$ with maps $\phi_{\alpha}:U_{\alpha}\longrightarrow V_{\alpha}\times W_{\alpha}$ where $V_{\alpha}$ is open in ${\mathbb{H}}^p:=\{(x_1,...,x_p)\in \R^p:x_1\geq 0\}$ and $W^{q}$ is open in $\R^q$. Changes of coordinates are in the form
 \begin{equation}v'=\phi(v,w),\quad w'=\psi(w)\end{equation} ($\psi$ is a local diffeomorphism).
 Such an atlas is assumed to be maximal among all collections of this type. 
 The integer $p$ is the dimension of the foliation, $q$ its codimension and $p+q=\operatorname{dim}(X_0).$
 In each foliated chart, the connected components of subsets as $\phi_{\alpha}^{-1}(V_{\alpha}\times \{w\})$ are called \underline{plaques}. The plaques coalesce to give maximal connected injectively immersed (not embedded) submanifolds called leaves. One uses the notation $\mathcal{F}_0$ for the set of leaves. Note that in general each leaf passes infinitely times trough a foliated chart so a foliation is only locally a fibration.
Taking the tangent spaces to the leaves one gets an integrable subbundle 

\noindent $T\mathcal{F}_0\subset TX_0$ that's transverse to the boundary i.e $T\partial X_0 + T\mathcal{F}_0=TX_0$ in other words the boundary is a submanifold transverse to the foliation.
\noindent Let given on $X$ (the manifold with cylinder attached) a smooth oriented foliation $\mathcal{F}$ with leaves of dimension $2p$ respecting the cylindrical structure i.e. 
\begin{enumerate}
\item The submanifold $\partial X_0$ is transversal to the foliation and inherits a $(2p-1,q)$ foliation $\mathcal{F}_{\partial}=\mathcal{F}_{|\partial X_0}$ with foliated atlas given by 

\noindent $\phi_{\alpha}:U_{\alpha}\cap \partial X_0 \longrightarrow \partial V_{\alpha} \times W_{\alpha}$. Note that the codimension is the same.
\item The restriction of the foliation on the cylinder is product type 

\noindent $\mathcal{F}_{|Z}=\mathcal{F}_{\partial}\times [0,\infty).$\end{enumerate}
\noindent These conditions imply that the foliation is normal to the boundary. We are going to introduce the notation for general Dirac type operators. We will specialize to the signature operator in the next section.
\noindent The orientation we choose is the one given by $(e_1,..,e_{2p-1},\partial_r)$ if $(e_1,..,e_{2p-1})$ is a positive leafwise frame for the induced boundary foliation. This is a way to specify the boundary Dirac type operator.
\noindent Let $E\longrightarrow X$ be a leafwise Clifford bundle with leafwise Clifford connection $\nabla^E$ and Hermitian metric $h^E$. Suppose each geometric structure is of product type on the cylinder meaning that if $\rho:\partial X_0\times [0,\infty) \longrightarrow \partial X_0$ is the base projection
$$E_{|Z}\simeq \rho^*(E_{|\partial X_0}),\quad h^E_{|\partial X_0}=\rho^*(h^{E}_{|\partial X_0}),\quad \nabla^E_{|Z}=\rho^*(\nabla^{E}_{|\partial X_0}).$$
\noindent Each geometric object restricts to the leaves to give a longitudinal Clifford module that's canonically $\mathbb{Z}_2$ graded by the leafwise chirality element. One can check immediately that the positive and negative boundary eigenbundles $E^+_{\partial X_0}$ and $E^-_{\partial X_0}$ are both modules for the Clifford structure of the boundary foliation.
\noindent Leafwise Clifford multiplication by $\partial_r$ induces an isomorphism of leafwise Clifford modules $c(\partial_r):E_{\partial X_0}^+\longrightarrow E_{\partial X_0}^-.$ Put $F=E^+_{|\partial X_0}$, the whole Clifford module on the cylinder $E_{|Z}$ can be identified with the pullback $\rho^*(F\oplus F)$ under the action:
tangent vectors to the boundary foliation $v\in T\mathcal{F}_{\partial}$ acts as
 $c^{E}(v)\simeq c^F(v)\Omega$ with $\Omega=\left(\begin{array}{cc}0 & 1 \\1 & 0\end{array}\right)$ while in the cylindrical direction $c^F(\partial_r)\simeq \left(\begin{array}{cc}0 & -1 \\1 & 0\end{array}\right)$. 
\noindent In particular one can form the longitudinal Dirac operator \footnote{the corresponding formula with $-\partial_r$, the inward pointing normal, is written to help the comparison with the orientation of A.P.S}
\begin{equation}
\label{diracco}
D=c(\partial_r)\partial_r+c_{|\mathcal F_0}\nabla^{E_{|\mathcal{F}_{\partial}}}=
c(\partial_r)\partial_r+\Omega D^{\mathcal{F}_{\partial}}=
c(-\partial_r)[-\partial_r-c(-\partial_r)\Omega D^{\mathcal{F}_{\partial}}].\end{equation}
Here 
$D^{\mathcal{F}_{\partial}}$ is the leafwise Dirac operator on the boundary foliation.  
\noindent In the following, these identifications will be omitted letting $D$ act directly on $F\oplus F$ according to
\begin{align}\nonumber
&\left(\begin{array}{cc}0 & D^- \\D^+ & 0\end{array}\right)=\left(\begin{array}{cc}0 & -\partial_r+D^{\mathcal{F}_{\partial}} \\ 
 \partial_r+D^{\mathcal{F}_{\partial}}
& 0\end{array}\right)=\left(\begin{array}{cc}0 & \partial_u+D^{\mathcal{F}_{\partial}} \\ 
 -\partial_u+D^{\mathcal{F}_{\partial}}
& 0\end{array}\right)
\end{align} where $u=-r$, $\partial_u=-\partial_r$ (interior unit normal). Remember that the signature operator is the Dirac operator corresponding to the natural Clifford module structure on the bundle of exterior algebras. We shall enter in details in the next section.

\noindent We assume the manifold is endowed with a holonomy invariant transverse measure $\Lambda$. Call $\mathcal{R}_0$ and $\mathcal{R}$ the equivalence relations of the foliations on $X_0$ and $X$ respectively both seen as measured groupoids with their natural Borel structure.

\section{The Hirzebruch formula}
The reference for the notation about the signature operator is the book by Berline Getzler and Vergne \cite{BeGeVe}.
Let $X$ be an oriented Riemannian manifold
and $|dvol|$ the unique volume form compatible with the metric i.e. the one assuming the value 1 on each positive oriented orthonormal frame. In other words $|\operatorname{dvol}|=|\sqrt g dx|.$
Define the Hodge $\ast$ operator in the usual way, $\ast e^{i_1}\wedge\cdot \cdot \cdot \wedge e^{i_k}= \operatorname{sign}({\sigma})e_{j_1}\wedge\cdot \cdot \cdot \wedge e_{i_{n-k}}$
where $(e_1,...,e_n)$ is an oriented orthonormal basis, $(i_1,...,i_k)$ and $(j_i,...,j_k)$ are complementary multindices and $\sigma$ is the permutation $\sigma:=\left(\begin{array}{cccccc}1 & . & . & . & . & n \\i_1 & . & i_k & j_1 & . & j_{n-k}\end{array}\right).$
 Since $\ast^2=(-1)^{|\cdot|(n-|\cdot|)}$ this is an involution on even dimensional manifolds.

\noindent The bundle $\Lambda{T^{*}X}$ of exterior algebras of $X$ is a natural Clifford module under the action defined by
$\label{cliffaction}c(e^{i}):=\epsilon(e_i)-\iota(e^{i})$ where  $\epsilon(e^i)\omega=e^{i}\wedge \omega$ is the exterior multiplication by $e^{i}$ and $\iota(e_i)$ is the contraction defined by the tangent vector $e_i$. These are related by the metric adjunction, $\epsilon(e^i)^*=\iota(e_i)$. The chirality involution is $\tau:=i^{[(n+1)/2]}c(e_1)\cdot\cdot\cdot c(e_n)$ and is related to the Hodge duality operator by 
$\tau=i^{[(n+1)/2]}\ast (-1)^{n|\cdot|+\frac{|\cdot|(|\cdot|-1)}{2}},$
following from the identity (same deegree forms)
$$\int_X \alpha \wedge \tau\beta=(-1)^{n|\cdot|+|\cdot|(|\cdot|-1)/2}i^{[2n+1]/2}\int_{X}(\alpha,\beta)|dx|$$ while $\int_X \alpha \wedge \ast \beta=\int_X (\alpha,\beta)|dx|.$
 As a consequence one can write the adjoint of $d$ in two different ways,
$$d^*=-\ast d \ast (-1)^{n|\cdot|+n}=-(-1)^n \tau d \tau.$$ Sections of the positive and negative eigenbundles of $\tau$ are called the \underline{self--dual} and \underline{anti self--dual} differential forms respectively and denoted by $\Omega^{\pm}(X).$

\noindent Now suppose $n$ is even, and $X$ is compact. The bilinear form on the middle cohomology $H^{n/2}(X;\R)$ defined by $(\alpha,\beta)\longmapsto \int_X \alpha \wedge \beta$ satisfies the identity 

\noindent $(\alpha,\beta)=(-1)^{n/2}(\beta,\alpha).$
In particular if $n$ is divisible by four this is \underline{symmetric} and has a signature $\sigma(X)$ i.e. the number $p-q$ entering in the representation
$Q(x)=x_1^2+\cdot \cdot \cdot + x_p^2- x_{p+1}^2-\cdot \cdot \cdot -x_q^2$
 of the associated quadratic form (this is independent by the choosen basis). In this situation the chiral Dirac operator $d+d^*$ acting on the space of differential forms is called the \underline{Signature operator}\footnote{it differs from the Gauss--Bonnet operator $d+d^*$ only for the choice of the involution}
 $$(d+d^*)=D^{\operatorname{sign}}=\left(\begin{array}{ccc}0 & D^{\operatorname{sign},-}  \\D^{\operatorname{sign},+} & 0  \end{array}\right):\Omega^+(X)\oplus \Omega^-(X)\longrightarrow \Omega^+(X)\oplus \Omega^-(X)$$
 \noindent The Atiyah--Singer index theorem specializes, for the signature, to the Hirzebruch formula
 $$\operatorname{ind}(D^{\operatorname{sign},+})=\sigma(X)=\int_X \operatorname{L}(X)$$ where $\operatorname{L}(X)$ is the $\operatorname{L}$--genus, $\operatorname{L}(X)=(\pi i)^{-n/2}\operatorname{det}^{1/2}\Big{(}\dfrac{R}{\operatorname{tanh(R/2)}}\Big{)}$ with $R$ the Riemannian curvature form.
The relationship between the Hirzebruch formula that admits a purely topological proof (based on cobordism) and the Atiyah Singer formula is given by the Hodge theorem stating a natural isomorphism between the space of {harmonic forms} $\mathcal{H}^q(X)$ i.e. the kernel of the forms laplacian $\Delta=(d+d^*)^2$ and the cohomology $H^q(X)$ together with Poincar\'e duality.
\bigskip
 
 \noindent Now on a $4k$--dimensional manifold with boundary $X_0$ with product structure the situation is much more complicated. The signature formula is the most important application of the index theorem in the A.P.S. paper. The operator can be written on a collar around the boundary as 
 $D^{\operatorname{sign},+}
 =\sigma(\partial_u+B)$ with the isomorphism $\sigma:\Omega(\partial X_0)\longrightarrow \Omega^+(X_0)$ and $B$ is the self--adjoint operator on $\Omega(\partial X_0)$ defined by $B\alpha=(-1)^{k+p+1}(\ast_{\partial}d-d\ast_{\partial})\alpha$ with $\epsilon(\alpha)=\pm 1$ according to $\alpha$ even or odd degree while $\ast_{\partial}$ is the Hodge duality operator on $\partial X_0$. Since $B$ commutes with
 
\noindent $\alpha \mapsto (-1)^{|\alpha|}\ast_{\partial}\alpha$ and preserves the parity of forms we have the splitting $B=B^{\textrm{ev}}\oplus B^{\textrm{odd}}$ and the dimension of the kernel at the boundary as the $\eta$ invariant are twice that of $B^{\textrm{ev}}$. The A.P.S index theorem says
 $$\operatorname{ind}(D^{\operatorname{sign},+}
)=\underbrace{h^+-h^-}_{\operatorname{ind}_{L^2}(D^{\operatorname{sign},+})}-h^-_{\ty}=
\int_XL-h(B^{\operatorname{ev}})-\eta(B^{\operatorname{ev}})$$  \noindent where $h^{\pm}$ are the dimensions of the $L^2$--harmonic forms on the manifold ${X}$ with a cylinder attached and $h^{-}_{\ty}$ is the dimension of the limiting values of the extended $L^2$--harmonic forms in $\Omega^-(X)$.

\noindent The identifications of all these numbers with topological quantities require some work. 
\begin{enumerate}
\item The space $\mathcal{H}({X})$ of $L^2$--harmonic forms on the elonged manifold ${X}$ is naturally isomorphic to the image $\widehat{H}(X_0)$ of the natural map\footnote{the inclusion of the compact support cohomology into the ordinary one}

\noindent $H^{*}_{0}({X})\longrightarrow H^{*}({X}).$ Equivalently  one can use the relative de Rham cohomology
$H^*(X_0,\partial X_0)\longrightarrow H^*(X_0)$
defined with boundary conditions $\omega_{|\partial X_0}=0$ on the de Rham complex. This is the role played by Hodge theory in the boundary case.
\item The signature $\sigma(X_0)$ of a manifold with boundary is defined to be the signature of the non--degenerate quadratic form on the vector space
$\widehat{H}^{2k}(X_0).$ This is induced by the degenerate quadratic form given by the cup--product on the relative cohomology 
$H^{2k}(X_0,\partial X_0)$. By Lefschetz duality the radical of this quadratic form is exactly the kernel of the mapping 

\noindent $H^{2k}(X_0,\partial X_0)\longrightarrow H^{2k}(X_0)$ then
$\sigma(X_0)=h^+-h^-=\operatorname{ind}_{L^2}(D^{\operatorname{sign},+}).$
\item Finally A.P.S get rid of the third number $h^-_{\ty}$ proving that 

\noindent $h^-_{\ty}=h^+_{\ty}=h(B^{\textrm{ev}})$ that together with $h^{+}_{\ty}+h^-_{\ty}=2h(B^{\textrm{ev}})$ gives the final signature formula
\\ \noindent $\sigma(X_0)=\int_{X_0}L-\eta(B^{\operatorname{ev}}).$

\end{enumerate}
\section{Computations with the leafwise signature operator}\label{computation}
\noindent Let start with a compact manifold with boundary $X_0$ equipped with an oriented $4k$--dimensional foliation transverse to the boundary. Suppose every geometric structure to be of product type near the boundary.
As usual attach an infinite cylinder $Z_0=\partial X_0\times [0,\ty)_r$ and extend all the geometric datas.
 The leafwise signature operator corresponds to the leafwise Clifford action defined above on the leafwise exterior bundle $\Lambda T^*{\mathcal{F}}$. If $(e_1,...,e_{4k-1},\partial_r)$ is a leafwise positive orthonormal frame near the boundary, the leafwise chirality element \footnote{we omit simbols denoting leafwise action for ease of reading} satisfies
\begin{align*} 
\tau:=i^{2k}c(e^1)\cdot \cdot \cdot c(e^{4k-1})c(dr)&=i^{2k}\ast (-1)^{|\cdot|(|\cdot|-1)/2}\\&=-i^{2k}c(dr)c_{\partial}=-i^{2k}c(dr)\ast_{\partial}(-1)^{|\cdot|+|\cdot|(|\cdot|-1)/2}\end{align*}
where $\ast$ is leafwise Hodge duality operator, the element $c_{\partial}=c(e^1)\cdot \cdot \cdot c(e^{4k-1})$ is a part for the $i^{2k}$ factor the leafwise boundary chirality operator and $\ast_{\partial}$ is the leafwise boundary Hodge operator. 
\noindent On the cylinder the leafwise bundle $\Lambda T^*\mathcal{F}$ is isomorphic to the pulled back bundle $\rho^*(\Lambda T^*\mathcal{F}_{\partial})$ (the projection on the base $\rho$ will be omitted throughout) while separating the $dr$ component on leafwise forms $\alpha=\omega + \beta \wedge dr$ yields an isomorphism
\begin{equation}\label{isobundles}(\Lambda T^*\mathcal{F})_{|\partial X_0} \longrightarrow (\Lambda T^* \mathcal{F}_{\partial})\oplus (\Lambda T^* \mathcal{F}_{\partial}),\end{equation} sometimes we shall write $(\Lambda T^*\mathcal{F}_{\partial})\wedge dr$ for the second addendum in \eqref{isobundles} to remember this isomorphism. An easy computation involving rules as 

\noindent $d\omega =d_{\partial}\omega+(-1)^{|\omega|}\partial_r \omega\wedge dr$ for $\omega \in C^{\ty}([0,\ty);\Lambda T^*\mathcal{F}_{\partial})$ and 

\noindent $c(dr)(\omega+\alpha \wedge dr)=(-1)^{|\omega|}\omega \wedge dr-(-1)^{|\alpha|}\alpha$ shows that w.r.t. the direct sum \eqref{isobundles},
\begin{equation}\label{repp1}D^{\operatorname{sign}}=\left(\begin{array}{cc}d_{\partial}+c_{\partial}d_{\partial}c_{\partial} & -(-1)^{|\cdot|}\partial_r \\(-1)^{|\cdot|} \partial_r & c_{\partial}d_{\partial}c_{\partial}\end{array}\right)=c(dr)\partial_r+(d_{\partial}+c_{\partial}d_{\partial}c_{\partial})\oplus (d_{\partial}+c_{\partial}d_{\partial}c_{\partial})
\end{equation} and 
\begin{equation}\label{tauchib}\tau=i^{2k}\left(\begin{array}{cc}0 & c_{\partial}(-1)^{|\cdot|} \\
-c_{\partial}(-1)^{|\cdot|}
 & 0\end{array}\right)
.\end{equation}
Since $d_{\partial}^*=\tau_{\partial}d_{\partial}\tau_{\partial}=c_{\partial}d_{\partial}c_{\partial}$ formula \eqref{repp1} is equivalent to
$D^{\operatorname{sign}}=c(dr)\partial_r+(d_{\partial}+d_{\partial}^*)\oplus(d_{\partial}+d_{\partial}^*).$
\noindent There's also another important formula corresponding to the fact that $d+d^*$ anticommutes with $\tau$. Denote $\Omega^{\pm}(\mathcal{F})$ the positive (negative) eigenbundles i.e. the bundles of leafwise auto--dual (anti auto--dual) forms. We can write the operator on the cylinder as an operator on sections of the direct sum $\rho^*(\Omega^+(\mathcal{F})_{|\partial X_0}\oplus \Omega^+(\mathcal{F})_{|\partial X_0} )$ as the matrix
\begin{align}\nonumber (-1)^{|\cdot|}\partial_r\left(\begin{array}{cc}0 & -1 \\1 & 0\end{array}\right)+(\ast_{\partial}d_{\partial}-d_{\partial}\ast_{\partial})i^{2k}(-1)^{|\cdot|}(|\cdot|-1)/2\left(\begin{array}{cc}0 & 1 \\1 & 0\end{array}\right) 
\\
\nonumber
=
c(dr)\partial_r+(\ast_{\partial} d_{\partial}-d_{\partial}\ast_{\partial})i^{2k} (-1)^{|\cdot|(|\cdot|-1)/2}\Omega.
\end{align}
\noindent To pass from one representation to another we have to consider the following compositions
\begin{equation*}
\xymatrix{ \Lambda T^*{\mathcal{F}_\partial} \ar[r]^-{i_1} & (\Lambda T^*{\mathcal{F}_\partial})\bigoplus (\Lambda T^*{\mathcal{F}_\partial})\wedge dr
\ar[r]^-{  1+\tau}
&
\Omega^+(\mathcal{F})
 \ar[r]^{d+d^*}&
 \Omega^-(\mathcal{F})\ar[r]^-{\operatorname{Pr}_2}&
\Lambda T^*{\mathcal{F}_\partial}}.
\end{equation*}
\begin{equation*}
\xymatrix{ \Lambda T^*{\mathcal{F}_\partial} \ar[r]^-{i_2} & \Lambda (T^*{\mathcal{F}_\partial})\bigoplus (\Lambda T^*{\mathcal{F}_\partial})\wedge dr
\ar[r]^-{  1-\tau}
&
\Omega^-(\mathcal{F})
 \ar[r]^-{d+d^*}&
 \Omega^+(\mathcal{F})\ar[r]^-{\operatorname{Pr}_1}&
\Lambda T^*{\mathcal{F}_\partial}.
}\end{equation*}
where $i_j$ is the inclusion on the $j$--th factor and $\operatorname{Pr}_{j}$ is the corresponding projection.

\section{The Analytic signature}
\noindent We can immediately define the analytic signature. It is the $L^2$ measured (chiral) index of the signature operator on the elongated manifold. In the paper \cite{io} we proved this is well defined and finite. In particular the $\Lambda$--dimensions of the extended spaces of solutions are finite.
\begin{dfn}
The $\Lambda$--analytic signature of the foliated manifold with boundary $X_0$ is the measured $L^2$--chiral index of the signature operator on the foliated manifold with a cylinder attached
\begin{equation}\label{2111}
\sigma_{\Lambda,\operatorname{an}}(X_0,\partial X_0):=\operatorname{ind}_{L^2,\Lambda}(D_X^{\operatorname{sign},+}).\end{equation}
Only here we write $(D_X^{\operatorname{sign},+})$ to stress that we consider the leafwise signature operator on $X$.
\end{dfn}
\noindent From the A.P.S index formula proved in \cite{io} and the standard identification of the Atiyah--Singer integrand for the signature operator \cite{BeGeVe}, formula \eqref{2111} becomes
$$\sigma_{\Lambda,\operatorname{an}}(X_0,\partial X_0)=\langle L(X),C_{\Lambda}\rangle +1/2[\eta_{\Lambda}(D^{\mathcal{F}_{\partial}})-h^+_{\Lambda}+h^{-}_{\Lambda}]$$ where $L(X)$ is the tangential $L$--characteristic class, the numbers $h^{\pm}_{\Lambda}$ and the foliation eta--invariant are referred to the boundary operator.

\noindent As in \cite{AtPaSi1} first we have to identify these numbers. Minor modifications of the proof of Vaillant \cite{Vai} (a complete proof  in \cite{iosuarxiv}) are needed in order to prove that for the signature operator
$h^+_{\Lambda}=h^-_{\Lambda}.$ Consequently the  signature formula  reduces to 
\begin{equation}\label{ffo}\sigma_{\Lambda,\,\operatorname{an}}(X_0,\partial X_0)=\langle L(X),C_{\Lambda}\rangle +1/2[\eta_{\Lambda}(D^{\mathcal{F}_{\partial}})].\end{equation}

\section{Fields of sesquilinear forms}
We shall give some definitions about Borel fields of sesquilinear forms in the setting of non commutative integration theory \cite{Cos}.
Let $\xymatrix{{\mathcal{G}}\ar^r@<2pt>[r]\ar_s@<-2pt>[r] &\mathcal{G}^{(0)}}$ be a Borel groupoid with a square integrable representation on a Borel field of Hilbert spaces $(H_x)_{x\in \mathcal{G}^{(0)}}$. In the next it will be $\mathcal{G}=\mathcal{R}_0$ or $\mathcal{G}=\mathcal{R}.$ Let $q=(q_x)_{x\in \mathcal{G}^{(0)}}$ be a $\mathcal{G}$--equivariant field of sesquilinear symmetric forms\footnote{a sesquilinear form $q$, antilinear in the second variable is symmetric if $\overline{q(\xi,\eta)}=q(\eta,\xi)$}, $q_x:H_x\times H_x\longrightarrow \mathbb{C}.$ By the Riesz lemma there exist a family of bounded selfadjoint intertwining operators $B=(B_x)_{x\in \mathcal{G}^{(0)}}$ such that $q_x(\xi,\eta)=(\xi,B_x \eta)_x$ for every $x\in \mathcal{G}^{(0)}.$ Measurability properties of $q$ are addressed, by definition to that of $B$. Now $B$ determines a field of hortonormal splittings $H_x=V_x^+\oplus V_x^0 \oplus V_x^-$ where $V_x^+$ $(V_x^-)$ is the image of the spectral projection $\chi_{(0,\t)}(B_x)$ ($\chi_{(-\ty,0)}(B_x)$) and $V_x^0$ is the kernel of $B_x$. If $\Lambda$ is a transverse measure on $\mathcal{G}$ one can measure the Random Hilbert spaces $V_x^+$ and $V_x^-$. If one of these formal dimension is finite one can define the $\Lambda$--\underline{signature} of $q$ to be $$\operatorname{sign}_{\Lambda}(q):=\operatorname{dim}_{\Lambda}(V^+)-\operatorname{dim}_{\Lambda}(V^-).$$
\section{The Hodge signature}
\noindent Consider the field of Hilbert spaces of $2k$ square integrable harmonic forms
\\\noindent $x\longmapsto \mathcal{H}_x^{2k}:= \operatorname{ker}\{\Delta^q_x:L^2(\Lambda^{2k} T^*L_x)\longrightarrow L^2(\Lambda^{2k} T^*L_x)\}$
where $L_x$ is a leaf of the $4k$--dimensional oriented foliation on the manifold with cylindrical ends $X$. Since leafwise harmonic forms are closed this is a field of subspaces of the field of the de Rham cohomologies $H^*(L_x)$. It inherits the structure of a measurable field of Hilbert spaces. This defines a square integrable representation of $\mathcal{R}$ and there is a natural field of symmetric forms $s^{\ty}_x:\mathcal{H}^{2k}_x \times \mathcal{H}^{2k}_x\longrightarrow \mathbb{C}$ defined by $$s^{\ty}_x(\alpha,\beta):=\int_{L_x}\alpha \wedge \beta= \int_{L_x}(\alpha,\ast \beta).$$ In the paper \cite{io} is proven that its $\Lambda$--signature is well defined. In fact this is precisely a statement about the finite $\Lambda$--dimensionality of the kernel of the leafwise Laplace Beltrami operator on $X$.

\begin{dfn}
The signature on harmonic forms (The Hodge signature or the harmonic signature) on the foliated elongated manifold is 
$$\sigma_{\Lambda,\operatorname{Hodge}}(X_0,\partial X_0):=\operatorname{dim}_{\Lambda}V^+-\operatorname{dim}_{\Lambda}V^-.$$\end{dfn}

\section{{Analytical signature}={Hodge signature}}
We prove that the analytical signature and the Hodge signature are equal. The boundary operator here is $B=\ast_{|\Omega^{2k}}$. Since the dimension of the foliation is $4k$ we have $\tau_{|\Omega^{2k}}=\ast_{|\Omega^{2k}}$. It follows that $V^{\pm}=\operatorname{ker}_{L^2}(D^{\operatorname{sign},\pm}).$ Then thanks to the index formula in the paper \cite{io} we have this first result
\begin{theorem}
The (measured) analytical signature and the (measured) harmonic signature coincide,
\begin{equation}\label{anasign}\sigma_{\Lambda,\operatorname{an}}(X_0,\partial X_0)=\sigma_{\Lambda,\operatorname{Hodge}}(X_0,\partial X_0)=\langle L(X),C_{\Lambda}\rangle +1/2[\eta_{\Lambda}(D^{\mathcal{F}_{\partial}})].\end{equation}
\end{theorem}

\section{The $L^2$--de Rham signature} The goal of this section is to give the definition of the de Rham signature for the foliated manifold with boundary and the proof that this is equal to the harmonic signature.
\subsubsection{Manifolds with boundary with bounded geometry}\label{bbg}
\noindent The generic leaf of $(X_0,\mathcal{F})$ is a Riemannian manifold with boundary with bounded geometry \cite{scht,sch1,sch2}.
\begin{dfn}\label{deltabounded}
We say that a manifold with boundary $M$ equipped with a Riemannian metric has bounded geometry if the following properties hold true.
\begin{description}
\item[Normal collar]: there exists $r_C>0$ so that the geodesic collar

\noindent $N:=[0,r_{C})\times \partial M:(t,x)\longmapsto \operatorname{exp}_x(t\nu_x)$
is a diffeomorphism onto its image. Here $\nu_x$ is the unit inward normal vector at $x\in \partial M$. Equip $N$ with the induced metric. In the sequel $N$ and its image will be identified. Denote $\operatorname{im}[0,r_C/3)\times \partial M$ by $N_{1/3}$. \item[Injectivity radius of the boundary]: the injectivity radius of $\partial M$ is positive, $r_{\operatorname{inj}}(\partial M)>0$
\item[Injectivity radius of $M$]: there is $r_i>0$ so that for $x\in M- N_{1/3}$ the exponential mapping is a diffeomorphism on $B(0,r_1)\subset T_xM$. In particular if we identify $T_xM$ with $\R^m$ via an orthonormal frame we have Gaussian coordinates $\R^m \supset B(0,r_i)\longrightarrow M$ around any point in $M- N_{1/3}$
\item[Curvature bounds]: for every $K\in \mathbb{N}$ there is some $C_K>0$ so that $|\nabla^i R|\leq C_K$ and $|\nabla^{\partial}l|\leq C_K$, $0\leq i\leq K$. Here $\nabla$ is the Levi--Civita connection on $M$, $\nabla^{\partial}$ is the Levi--Civita connection on $\partial M$ and $l$ is the second fundamental form tensor with respect to $\nu$.
\end{description}
\end{dfn}
\noindent Choose some $0<r_1^C<r_{\operatorname{inj}}(\partial M)$. Near points $x' \in \partial M$ on the boundary one can define {\bf{normal collar coordinates}} by iteration of the exponential mapping of $\partial M$ and that of $M$,
$$k_{x'}:\underbrace{B(0,r_i^C)}_{\subset \R^{m-1}}\times [0,r_C)\longrightarrow M, (v,t)\longmapsto \operatorname{exp}^{M}_{\operatorname{exp}_{x'}^{\partial M}(v)}(t\nu).$$ For points $x\in M-N_{1/3}$ standard {\bf{Gaussian coordinates}} are defined via the exponential mapping. In the following we shall call both {\bf{normal coordinates}}. It is a non trivial fact that the condition on curvature bounds in definition \ref{deltabounded} can be substituted by uniform control of each derivative of the metric tensor $g_{ij}$ and its inverse $g^{ij}$ on normal coordinates. 

\noindent The definition extends to bounded geometry vector bundles on boundary manifolds with bounded geometry and each object of uniform analysis like i.e. uniformly bounded differential operators \cite{sch2}.
In particular, using a suitable partition of the unity adapted to normal coordinates one can define uniform Sobolev spaces (different coordinates give equivalent norms so we get hilbertable spaces) and every basic result continues to hold.
\begin{prop}\label{schick}
Let $E\longrightarrow M$ a bundle of bounded geometry over $M$. Suppose $F \longrightarrow \partial M$ is a bounded vector bundle. 
We have the following properties for the uniform Sobolev spaces $H^s(E), H^t(F)$ of sections ($s,t\in \R$).
\begin{enumerate}
\item $H^s(E), H^t(F)$ is an Hilbert space (inner product depending on the choices).
\item The usual (bounded) Sobolev embedding theorem holds with values on the Banach space $C_b^k(E)$ of all sections with the first $k$ derivatives uniformly bounded in normal coordinates.
$$H^s(E)\hookrightarrow C^k_b(E),\quad \mbox{whenever} \quad s>m/2+k.$$
\item For the bundle of differential forms one can use as Sobolev norm the one coming from the integral of the norm of covariant differentials

\noindent $\|\omega\|_k^2:=\sum_{i=0}^{k}\int_M\|\nabla^i \omega(x)\|^2_{T^{*}_x M \otimes \Lambda T^*M}|dx|.$
\item For $s<t$ we have a bounded embedding with dense image $H^t(E)\subset H^s(E)$. The map is compact if and only if $M$ is compact. We can define, as usual
$$H^{\ty}(E):=\bigcap_s H^s(E),\quad H^{-\ty}(E):=\bigcup_s H^s(E).$$
\item Let $p:C^{\ty}(E)\longrightarrow C^{\ty}(F)$ a $k$--bounded boundary differential operator i.e the composition of an order $k$ bounded differential operator on $E$ with the morphism of restriction to the boundary. Then $p$ extends to be a bounded operator
$p:H^s(E)\longrightarrow H^{s-k-1/2}(F)$ if $s>k+1/2.$
\noindent In particular we have the bounded restriction map $H^s(E)\longrightarrow H^{s-1/2}(E_{|\partial M})$ if $s>1/2$.
\item $H^s(E)$ and $H^{-s}(E)$ are dual to each according to the  extension by continuity of the pairing
$$(f,g)=\int_Mg(f(x))|dx|; \, f\in C^{\ty}_0(E), \,g\in C^{\ty}_0(E^*)$$ where $E^*$ is the dual bundle of $E$.
\noindent
If $E$ is a bounded Hermitian or Riemannian bundle, then the norm on $L^2(E)$ defined by charts is equivalent to the usual $L^2$--norm
$$|f|^2:=\int_M(f,f)_x|dx|,\,f\in C^{\ty}_0(E).$$
Moreover $H^s(E)$ and $H^{-s}(E)$ are dual to each other by extension of 

\noindent $(f,g)=\int_M(f,g)_x|dx|.$
\end{enumerate}
\end{prop}
\subsubsection{Random Hilbert complexes}\label{derham}
\noindent We are going to define the $L^2$--de Rham complexes along the leaves. These are particular examples of Hilbert complexes studied in complete generality in \cite{brun}. So let $x\in X_0$, consider the unbounded operator with Dirichlet boundary conditions $$d_{L_x^0}:\Omega^k_{d,x}=\{\omega \in C^{\ty}_0(\Lambda T^kL_x^0); \omega_{|\partial L_x^0}=0\}\subset L^2_x(\Lambda T^kL_x^0)\longrightarrow L^2_x(\Lambda T^kL_x^0).$$
\noindent Being a differential operator it is closable. Let $A^k_x(L_x^0,\partial L_x^0)$ the domain of its closure i.e the set of $L^2$ limits $\omega$ of sequences $\omega_n$ such that also the $d\omega_n$ converges in $L^2$ to some $\eta=:d\omega$. 
The graph norm $\|\cdot\|_A^2:=\|\cdot \|_{L^2}^2+\|d\cdot \|_{L^2}^2$ gives the graph the structure of an Hilbert space where $d$ is bounded. 
It is easily checked that $d(A^k_x)\subset \operatorname{ker}(d:A^{k+1}_x\longrightarrow L^2_x)$.  In particular we have a Hilbert cochain complex
$\cdot\cdot\cdot \longrightarrow A^{k-1}_x\longrightarrow A^k_x \longrightarrow A^{k+1}_x\longrightarrow \cdot \cdot \cdot$
\noindent with
\begin{itemize}
\item Cycles $Z^k_x(L_x^0,\partial L_x^0):=\operatorname{ker}(d:A^k_x\longrightarrow A^{k+1}_x)$.
\item Boundaries \noindent $B^k_x(L_x^0,\partial L_x^0):=\operatorname{range}(d:A^{k-1}_x\longrightarrow A^{k}_x).$
\end{itemize}
\begin{dfn}
 The $L^2$ (reduced )\footnote{the word reduced stands for the fact we use the closure to make the quotient, also the non reduced cohomology can be defined. For a $\Gamma$ covering of a compact manifold the examination of the difference reduced/unreduced cohomology leads to the definition of the Novikov--Shubin invariants\cite{luk}} relative de Rham cohomology of the leaf $L_x^0$
 is defined by the quotients
 $$H^{k,x}_{dR,(2)}(L_x^0,\partial L_x^0):={Z^k_x(L_x^0,\partial L_x^0)}\Big{/}{\overline{B^k_x(L_x^0,\partial L_x^0)}}.$$
\end{dfn}
\noindent One take the closure in the definition in order to assure the quotient to be an Hilbert space. 
\noindent Similarly the $L^2$--de Rham cohomology of the whole leaf, $H^{k,x}_{dR,(2)}(L_x^0)$ is defined using no (Dirichlet) boundary conditions. In particular $A_x^k(L_x^0)$ will be used to denote the domain of the closure of the differential as unbounded operator on $L^2(L_x^0)$ defined on compactly supported sections (the support possibly meeting the boundary). The subscript $dR$ helps to make distinction with Sobolev spaces.
\noindent 
Each one of this spaces is naturally isomorphic via $L^2$--Hodge theory to a corresponding space of harmonic forms. 
\begin{dfn}
The space of degree $k$--$L^2$--
harmonic forms with Dirichlet boundary conditions on $\dl$ is 
$$\mathcal{H}^k_{(2)}(L_x^0,\dl):=\{\omega \in C^{\ty}\cap L^2,\, \omega_{|\dl}=0,\, (\delta \omega)_{|\dl}=0,\,{(d\omega)_{|\dl}=0} \}$$
\end{dfn}\noindent The condition $(d\omega)_{|\dl}=0$ can be omitted in the definition since is automatically satisfied. The boundary conditions are exactly the square of the Dirichlet boundary condition for the Dirac operator $d+\delta$.
\noindent Since each leaf is complete a generalization of an idea of Gromov shows that these forms are closed and co--closed \cite{scht,sch1}
$$\mathcal{H}^k_{(2)}(L_x^0,\dl)=\{\omega \in C^{\ty}\cap L^2(\Lambda^kL_x^0), d\omega=0,\,\delta \omega=0,\,\omega_{|\dl}=0\}.$$ 
\noindent Furthermore there's the $L^2$--orthogonal Hodge decomposition \cite{scht,sch1}
$$L^2(\Lambda^k T^*L_x^0)=\mathcal{H}^k_{(2)}(L_x^0,\dl)\oplus \overline{d^{k-1}\Omega_{d,x}^{k-1}(L_x^0,\dl)}^{L^2}\oplus \overline{\delta^{k+1}\Omega_{\delta,x}^{k+1}(L_x^0,\dl)}^{L^2} $$ 
where 
$\Omega_{d,x}^{k-1}:=\{\omega \in C^{\ty}_0(\Lambda^{k-1} T^*L_x^0),\,\omega_{|\dl}=0\}$
and the corresponding one for $\delta$ (with no boundary conditions) 
$\Omega_{\delta,x}^{k+1}:=\{\omega \in C^{\ty}_0(\Lambda^{k+1} T^*L_x^0)\}.$
 These decompositions show with a little work that the inclusion $\mathcal{H}^k(L_x^0,\dl)\hookrightarrow A^k_x$ induces isomorphism in cohomology (Hodge--de Rham Theorem)
$$\mathcal{H}^k(L_x^0,\dl)\cong H_{dR,(2)}^k(L_x^0,\dl).$$ This is a consequence of the fact that the graph norm (of $d$) and the $L^2$ norm coincide on the space of cycles $Z^k_x$.  Similar Hodge isomorphisms hold for the non--relative spaces.

\noindent As $x$ varies in $X_0$ they form measurable fields of Hilbert spaces. We discuss this aspect in a slightly more general way applicable to other situations. Remember that a measurable structure on a field of Hilbert spaces over $X_0$ is given by a fundamental sequence of sections, $(s_x)_{x\in X_0}$, $s_n(x)\in H_x$ such that $x\longmapsto \|s_n(x) \|_{H_x}$ is measurable and $\{s(x)\}_n$ is total in $H_x$ (see chapter IV in \cite{take} ). \begin{prop}If for a family of closed densely defined operators $(P_x)$ with minimal domain $\mathcal{D}(P_x)$ a fundamental sequence $s_n(x) \in \mathcal{D}(P_x)$ is a core for $P_x$ and $P_x s_n(x)$ is measurable for every $x$ and $n$ then the family $P_x$ is measurable in the sense of closed unbounded operators i.e. the family of projections $\Pi_x^g$ on the graph is measurable in the square field $\{H_x\oplus H_x\}_x$ with its product measurable structure.
\end{prop}
\begin{proof}\noindent Since the graph is generated by vectors $(s_n(x),P_x s_n(x))$ then the field of projections is measurable.
\end{proof}
\noindent The above lemma can be applied to the $(A^{k}_x(L_x^0,\dl))_x$ in fact in the appendix of \cite{hl} a fundamental sequence $\varphi_n$ of sections with the property that each $(\varphi_n(\cdot))_{|L_x^0}$ is smooth and compactly supported is showen to exist. The same proof works for manifold with boundary furthermore since the boundary has zero measure one can certainly require to each $\varphi_n$ to vanish on the boundary.

\noindent In particular we have defined complexes of square integrable representations. Reduction modulo $\Lambda$ almost everywhere gives complexes of random Hilbert spaces (with unbounded differentials) for which we introduce the following notations, 
\begin{itemize}
\item $(L^{2,\mathcal{F}}(\Omega^{\bullet}X_0),d)$ is the complex of Random Hilbert spaces obtained by reduction 
$\Lambda$--a.e.
from the field of Hilbert complexes
\begin{equation}\label{sssd}
\xymatrix{
{\cdot\cdot\cdot} \ar[r]
&L^2(\Lambda^k T^*L_x^0)\ar[r]^-d
& L^2(\Lambda^{k+1}T^* L_x^0)\ar[r]&\cdot \cdot \cdot
}
\end{equation}
\item $H^{\bullet,\mathcal{F}}_{dR,(2)}(X_0)$ is the random Hilbert space obtained by $\Lambda$--a.e. reduction from the reduced $L^2$--homology of
\eqref{sssd}.
\item $(L^{2,\mathcal{F}}(\Omega^{\bullet}X_0,\partial X_0),d)$ is the complex of Random Hilbert spaces obtained by $\Lambda$--a.e. reduction from the field of Hilbert complexes with Dirichlet boundary condition
\begin{equation}\label{ssssd}
\xymatrix{
{\cdot\cdot\cdot} \ar[r]
&L^2(\Lambda^k T^*L_x^0)\ar[r]^-d
& L^2(\Lambda^{k+1}T^* L_x^0)\ar[r]&\cdot \cdot \cdot
}
\end{equation}
with differentials considered as unbounded operators with domain $A_x^{k}(L_x^0,\partial L_x^0).$
\item $H^{\bullet,\mathcal{F}}_{dR,(2)}(X_0,\partial X_0)$ is the random Hilbert space of the reduced homology of \eqref{ssssd}.
\end{itemize}
\subsubsection{Definition of the $L^2$--de Rham signature}

\noindent Remember that $\operatorname{dim}(\mathcal{F})=4k$. Consider the measurable field of Hilbert spaces $A^k_{x}(L_x^0,\dl)$ of the minimal domains of the de Rham leafwise differential with Dirichlet boundary conditions $\omega_{| \dl}=0$ as above, with the graph Hilbert structure and the induce Borel structure. This square integrable representation of $\mathcal{R}_0$ carries a field of bounded symmetric sesquilinear forms defined by
$$s^0_x:A^{2k}_{x}(L_x^0,\dl)\times A^{2k}_{x}(L_x^0,\dl)\longrightarrow \mathbb{C},(\omega,\eta)\longmapsto \int_{L_x^0} \omega \wedge \overline{\eta}=\int_{L_x^0}(\omega,\ast \eta)d\nu^x.$$ This is the $\mathbb{C}$--antilinear (in the second variable) extension of the wedge product on forms. The complex conjugate is defined by $\overline{\sigma \otimes \gamma}=\sigma \otimes \bar{\gamma}$ and $\nu^x$ is the Leafwise Riemannian metric. Here also the scalar product $(\cdot,\cdot)$ on forms is extended to be sesquilinear.
\begin{lem}
The sesquilinear form $s^0_x$ passes to the $L^2$--relative cohomology of the leaf $H^{2k}_{dR,(2)}(L_x^0,\dl)$ and factorizes through the image of the natural map $H^{2k}_{dR,(2)}(L_x^0,\dl)\longrightarrow H^{2k}_{dR,(2)}(L_x^0)$ of the $L^2$--relative de Rham cohomology to the $L^2$--de Rham cohomology exactly as in the compact (one leaf) case.
\end{lem}
\begin{proof}
The first assertion is a consequence of the Stokes formula. Let 
$\omega \in A^{2k}_x(L_x^0,\dl)$ i.e. 
$\xymatrix{{\omega}_n\ar[r]^{L^2}&\omega}$, 
$\xymatrix{{d\omega}_n\ar[r]^{L^2}&0       }$ and 
$\theta_m\in C^{\ty}_0(\Lambda T^{2k-1}L_x^0)$, 
$\xymatrix{d\theta_m \ar[r]^{L^2}&\varphi}$ then 
$$s_x^0(\omega,\varphi)=\lim_{n,m}\int_{L_x^0}\omega_n \wedge \overline{d\theta_m}=\lim_{n,m}\int_{L_x^0}d(\omega_n \wedge \overline{\theta_m} )=\lim_{n,m}\int_{\dl} (\omega_n\wedge \theta_m)_{|\dl}=0.$$ The second one is clear and follows exactly from the classical case i.e. if 

\noindent $\beta_1=\beta_2+\lim_n d\rho_n$ with $\rho_n$ compactly supported with no boundary conditions write $$s_x^0([\alpha],[\beta])=s_x^0([\alpha],[\beta_2])+\lim_n \int \alpha \wedge \rho_n,$$ represent $\alpha$ as a $L^2$--limit of forms with Dirichlet boundary conditions and apply Stokes theorem again.
\end{proof}
\noindent For every $x$ the sesquilinear form $s_x^0$ on the cohomology corresponds to a bounded selfadjoint operator $B_x\in B(H_{dR,(2)}^{2k}(L_x^0,\dl     ))$ with $s_x^0(\alpha,\beta)=(\alpha, B_x \beta)$. Measurability properties of $(s_x^0)_{x\in X_0}$ are by definition measurability properties of the family $(B_x)_x.$ It is clear that everything varies in a Borel fashion (use again a smooth fundamental sequence of vector fields as in \cite{hl}). Then the $B_x$ define a self--adjoint random operator $B$ in the semifinite Von Neumann algebra $\operatorname{End}_{\Lambda}(H_{dR,(2)}^{2k}(X_0,\partial X_0))$ with trace $\operatorname{tr}_{\Lambda}$.

\begin{dfn}The $\Lambda$--$L^2$--de Rham signature of the foliated manifold with boundary $X_0$ is
$$\sigma_{\Lambda,dR}(X_0,\partial X_0):=\operatorname{tr}_{\Lambda}\chi_{(0,\ty)}(B)-\operatorname{tr}_{\Lambda}\chi_{(-\ty,0)}(B).$$
\end{dfn}

\section{{$L^2$--de Rham signature}$=${Hodge signature} }
\noindent The path to follow is clearly the one in the paper of of L{\"u}ck and Schick \cite{lus}. We shall show at the end of the section that we can reduce to the case in which
{\bf{every leaf meets the boundary}} or in other words {\bf{the boundary contains a complete transversal}}.
\subsubsection{The boundary foliation and $\mathcal{R}_0$ }\label{fundationtr}\label{primac} We refer to the Appendix \ref{appe} for a rapid account of Non commutative integration theory.
\noindent We have denoted by $\mathcal{F}_{\partial}$ the foliation induced on the boundary $\partial X_0$ i.e. the foliation whose a leaf is a connected component of the intersection of a leaf $L$ of $\mathcal{F}$ with the boundary. Let $\underline{\mathcal{R}_0}=\mathcal{R}(\mathcal{F}_{\partial})$ its equivalence relation with canonical inclusion $\underline{\mathcal{R}_0}\longrightarrow \mathcal{R}_0$. We are under the assumption that the boundary contains a complete transversal $T$. This is also a complete transversal for $\mathcal{F}_{\partial}$. Call $\nu_T$ its characteristic function \cite{Cos} on $\mathcal{R}_0$. Every transverse measure $\Lambda$ on $\mathcal{R}_0$ is univocally determined by the measure $\Lambda_{\nu_{T}}$ supported on $T$. As a consequence one gets a transverse measure, continue to call $\Lambda$, on $\underline{\mathcal{R}_0}$. \noindent Let now $(H,U)$ be a square integrable representation of $\mathcal{R}_0\longrightarrow X_0$ and $H$ its corresponding random Hilbert space, it pulls back to a square integrable representation $(H',U')$ of $\underline{\mathcal{R}_0}$. Also a random operator $A\in \operatorname{End}_{\Lambda}(H)$ defines by restriction a random operator $A'$ in $\operatorname{End}_{\Lambda}(H')$.
We are going to show that 
\begin{equation}\label{andrea}
\operatorname{tr}_{\Lambda}(A)=\operatorname{tr}_{\Lambda}(A').
\end{equation}
\noindent This is automatically proven if we think about the trace in terms of the operator valued weight $\int \operatorname{tr}_{H_x}(\cdot)d\Lambda_{\nu}(x)$ of Lemme 8 in \cite{Cos} (also Lemma \ref{11010} in the appendix). Of course we have to pay some care checking the domains of definitions of the two traces but from normality and square integrability the operators in form $\theta_{\nu}(\xi,\xi)$ as in Proposition 15 in \cite{Cos} (see also the appendix) furnish a sufficiently rich set to check the two. We prefer to look at the problem under a slightly different point of view. First remember that the trace of an operator is related to an integration of a Random variable (Proposition 14 Page 43 in \cite{Cos})  on $\mathcal{R}_0.$ So if one chooses as transverse function the characteristic function of $T$ and apply the Recipe of Connes finds out immediately that 
\begin{prop}
An intertwining operator between two square integrable representations of $R_0$ restricts to an intertwining operator between the pull--back representations of $\underline{\mathcal{R}_0}$ to give an element of the corresponding Von Neumann algebra with the same trace.
\end{prop}

\noindent This simple argument allows ourselves to consider, as an instrument short sequences  
$$\xymatrix{
0\ar[r]&A_x^{k-1}(L_x^0,\partial L_x^0)\ar[r]^-i&A_x^{k-1}(L_x^0)\ar[r]^-r &A_x^{k-1}(\partial L_x^0)\ar [r]&0,\,\,\,\,x\in \partial X_0
}$$
as sequences of Random Hilbert spaces associated to the boundary equivalence relation $\underline{\mathcal{R}_0}.$ In fact the third term seems not so naturally defined without passing to the boundary relation. 
It seems we have to say some words more about the relation between $\underline{\mathcal{R}_0}$ and ${\mathcal{R}_0}$ or, better its restriction to the boundary $(\mathcal{R}_0)_{|\partial X_0}$ . We shall investigate how $\underline{\mathcal{R}_0}$ sits inside $(\mathcal{R}_0)_{|\partial X_0}$ and the effect on the traces on the various algebras associated. Consider a class of $\underline{\mathcal{R}_0}$ i.e. a leaf of the boundary foliation; this is a connected component of a class of $(\mathcal{R}_0)_{|\partial X_0}$. In other words each class of $(\mathcal{R}_0)_{|\partial X_0}$ is a denumerable union of classes of $\underline{\mathcal{R}_0}$ i.e. the bigger one seems like to be some sort of denumerable union of the smaller under the obvious natural functor
$$\underline{\mathcal{R}_0}\longrightarrow (\mathcal{R}_0)_{|\partial X_0.}$$ In the measure theory realm denumerability means that $(\mathcal{R}_0)_{|\partial X_0}$ is not so bigger than $\underline{\mathcal{R}_0}$.
If one makes use of a complete transversal for $\underline{\mathcal{R}_0}$ to integrate natural\footnote{i.e. given by $L^2\bullet L$, where $L$ is left traslation on $\mathcal{R}_0$} Random Hilbert spaces associated to $(\mathcal{R}_0)_{|\partial X_0}$, this transversal touches denumerably times classes of $\underline{\mathcal{R}_0}$ hence  \emph{we are integrating (then taking traces) on the foliation induced on the boundary.} 
 The notion of properness for measurable functors helps to understand this intuitive fact. 
Recall from \cite{Cos} that a measurable functor 
 $F:\mathcal{G}\longrightarrow M$ with values standard measure spaces is proper if with respect to the diagram 
$$\xymatrix{
{\mathcal{G}}\ar[r]^-{F}\ar[d]^r&  X:=\bigcup_{x\in \mathcal{G}^{(0)}}\mathcal{G}^x\ar[dl]^\pi\\
\mathcal{G}^{(0)}}$$
$\mathcal{G}$ acts properly on $X$, i.e. there exist a strictly positive function $f\in\mathcal{F}^+(X)$ and a proper measure $\nu\in \mathcal{E}^+$ such that $\nu \ast f=1$. Here we recall the defining formula $$(\nu \ast f)(z):=\int_{\mathcal{G}^y} f(F({\gamma}^{-1})\cdot z)d\nu^y(\gamma).$$
Consider indeed the diagram
 $$\xymatrix{
{\underline{{\mathcal{R}_0}}}\ar[d]^{L'}\ar[r]^{L'}&X'=\bigcup_{x\in \partial X_0}(\mathcal{R}_0)^x_{|\partial X_0}\\
X=\bigcup_{x\in \partial X_0}\mathcal{R}_0^x}.$$ 
where $L'$ is the left multiplication functor $x\longmapsto (\mathcal{R}_0)_{|\partial X_0}^x$ while $L$ is left traslation in $\mathcal{R}_0$. Both are proper functurs because the first is the restriction of the multiplication of $\mathcal{R}_0$ and the second is the multiplication of the groupoid (Exemple after Definition 3 page 23 in \cite{Cos}). Associate to $L$ and $L'$  some local trace of an intertwining operator $B$ of a square integrable representation, say $x\longmapsto L^2(\partial L_x^0)$. We are saying that the target space $L'(x)$ is $(\mathcal{R}_0)^x_{|\partial X_0}$ and the measure is $f\longmapsto \alpha(f)=\operatorname{tr}(B^{1/2}fB^{1/2})$; the same association is done for $L$. Note that the integral $\int L' d\Lambda$ is exactly $\operatorname{tr}_{\Lambda}(B)$ in $\operatorname{End}_{\Lambda}(L^2(\partial L_x^0)).$ Now there is a Borel map associating to $z\in X$ a probability measure on $X'$ as in Proposition 4 pag 23 in \cite{Cos}.
It is the Dirac measure i.e. $z=(x,y)\longmapsto\lambda^z:=\delta_{(x,y)}.$ The property $$\int \lambda^zd\alpha^x(z)=\alpha'(x)$$ is, by definition immediately verified. Hence Proposition 4 pag. 23 in \cite{Cos} says that
\begin{prop}
The trace on endomorphisms of natural representations of $(\mathcal{R}_0)_{|\partial X_0}$ is equal to the trace arising from the foliation induced on the boundary.
\end{prop}
\subsubsection{Weakly exact sequences}
\noindent Consider for $x\in \partial X_0$ the Borel field of cochain complexes
$$\xymatrix{{}& {}\ar[d]^-d& {}\ar[d]^-d& {}\ar[d]^-d\\
0\ar[r]&A_x^{k-1}(L_x^0,\partial L_x^0)\ar[r]^-i\ar[d]^-d & A_x^{k-1}(L_x^0)\ar[r]^-r\ar[d]^-d & 
A_x^{k-1}(\partial L_x^0)\ar[r]\ar[d]^-d& 0\\
0\ar[r]& A_x^{k}(L_x^0,\partial L_x^0)\ar[r]^-i\ar[d]^d & A_x^{k}(L_x^0)\ar[r]^-r\ar[d]^-d & 
A_x^{k}(\partial L_x^0)\ar[r]\ar[d]^-d& 0\\
{}&{}&{}&{}}$$ where each morphism must be considered as an unbounded operator on the corresponding $L^2$. The inclusion $i$ is bounded since it is simply the restriction of the identity mapping on $L^2(L_x^0,\Lambda T^*L_x^0)$. The map $r$ is the restriction to the boundary.
\begin{prop}\label{kkkk}
\begin{enumerate}
\item For every $k$ the domain $A_x^{k}(L_x^0)$ is contained in the Sobolev space of forms $H^1(L_x^0,\Lambda T^*L_x^0)$. In particular the composition with $r$ makes sense.
\item The rows are exact. 
\end{enumerate}
\end{prop}
\begin{proof}
1. An element $\omega$ in $A_x^{k}(L_x^0)$ is an $L^2$--limit of smooth compactly supported forms $\omega_n$ with differential also convergent in $L^2$. Then since the Hodge $\ast$ is an isometry on $L^2$ also $\delta \omega_n=\pm \ast \omega \ast$ converges. In particular we can control the $L^2$--norm of $d\omega$ and $\delta \omega$. In particular we have control of the first covariant derivative. In fact $d+\delta=c \circ \nabla$ where $c$ is the (unitary) Clifford action. Then the second term can made less that the norm of $\nabla$ by bounded geometry.
In particular we have control on the order one Sobolev norm by proposition \ref{schick}. The remaining part follows from the fact that the restriction morphism is bounded from $H^1$ to $H^{1/2}\hookrightarrow L^2$. 

\noindent 2. The only non--trivial point is exactness in the middle but as a consequence of the bounded geometry the boundary condition on the first space extends to $H^1$ (see proposition $5.4$ in the thesis of Thomas Schick \cite{scht}) that together with point $1.$ gives exactness.
\end{proof}
\begin{rem}
Note that the proof of the above proposition also says that the induced morphisms $i_*$ and $r_*$ on the $L^2$--cohomology are bounded.
\end{rem}
\noindent Every arrow induces morphisms on the reduced $L^2$ cohomology. Miming the algebraic construction of the connecting morphism (everything works thanks to the above remark) we have, for every $x\in \partial X_0$ the long sequence of square integrable representations of the equivalence relation of the boundary foliation $\underline{\mathcal{R}_0}$
\begin{equation}\nonumber
\xymatrix{{\cdot \cdot \cdot}\ar[r]& H^{k,x}_{dR,(2)}(L_x^0,\partial L_x^0)\ar[r]^-{i_{*}} & 
H^{k,x}_{dR,(2)}(L_x^0)\ar[r]^-{r_{*}}&\\
{}&\ar[r]^-{r_{*}}&
H^{k,x}_{dR,(2)}(\partial L_x^0)\ar[r]^-{\delta}&
H^{k-1,x}_{dR,(2)}(L_x^0,\partial L_x^0)\ar[r]& \cdot \cdot \cdot}\end{equation}
\noindent Remove the dependence on $x$ to get a long sequence of Random Hilbert spaces over $\partial X_0$ with consistent notation with 
\eqref{sssd} and \eqref{ssssd}
\begin{equation}\label{success}
\xymatrix{{\cdot \cdot \cdot}\ar[r]& H^{k}_{dR,(2)}(X_0,\partial X_0)\ar[r]^-{i_{*}} & 
H^{k}_{dR,(2)}(X_0)\ar[r]^-{r_{*}}&\\
{}&\ar[r]^-{r_{*}}&
H^{k}_{dR,(2)}(\partial X_0)\ar[r]^-{\delta}&
H^{k-1}_{dR,(2)}(X_0,\partial X_0)\ar[r]& \cdot \cdot \cdot}
\end{equation}
\begin{dfn}\label{vonexact}
We say that a sequence of Random Hilbert spaces as \eqref{success} is $\Lambda$--weakly exact at some term if in the correspondig Von Neumann algebra of Endomorphisms the projection on the closure of the range of the incoming arrow coincides with the projection on the kernel of the starting one. As an example at point
$H_{dR,(2)}^k(X_0)$ we must have
$$\overline{\operatorname{range} i^*}=\operatorname{ker} i^*\in \operatorname{End}_{\Lambda}(H_{dR,(2)}^k(X_0)).$$
\end{dfn}
\noindent Of course such a sequence cannot be exact, just as in the case of Hilbert $\Gamma$--modules there are simple examples of non exacteness (see Example 1.19 in \cite{luk}, or the example on manifolds with cylindrical ends
 in the paper of Cheeger and Gromov \cite{ChGr}).
We shall see that a necessary condition to weakly exactness is (left) fredholmness, exactly as in the work of Cheeger and Gromov \cite{ChGr}.

\subsubsection{Spectral density functions and Fredholm complexes.}\label{1237}
\noindent Let $U$,$V$ two Random Hilbert spaces on $\mathcal{R}_0$ (for these facts also the holonomy groupoid or, more generally a Borel groupoid would work) and an {unbounded} 
{Random} {operator}
$f:\mathcal{D}(f)\subset U\longrightarrow V$ defined starting with a Borel family of closed densely defined operators $f_x:U_x\longrightarrow V_x$ intertwining the representation of $\mathcal{R}_0$. Since $f$ is closable, the question of measurability is addressed, as in \cite{nus} to the family of the closures. For every 
$\mu\geq 0$ let
$\mathcal{L}(f,\mu)$ be the set of measurable fields of subspaces $L_x\subset \mathcal{D}(f_x)\subset U_x$ such that for every $x\in X_0$ and $\phi\in L_x$, 
$
\|f_x(\phi)\|\leq \mu \|\phi\|$. After reduction modulo $\Lambda$ a.e. it becomes a set of Random Pre--Hilbert spaces we call 
$\mathcal{L}_{\Lambda}(f,\mu)$.
\begin{dfn}\label{frett}
The $\Lambda$--spectral density function of the family $\{f_x\}_x$ is the monotone increasing function
$\mu \longmapsto F_{\Lambda}(f,\mu):=\sup \{\operatorname{dim}_{\Lambda}:L\in \mathcal{L}_{\Lambda}(f,\mu)\}.$ Here of course one has to pass to the closure in order to apply the $\Lambda$--dimension. We say $f$ is \underline{$\Lambda$--Fredholm} if for some $\epsilon>0$, $F_{\Lambda}(f,\epsilon)<\ty$
\end{dfn}
\noindent We are going to show that this definition actually coincides with the (standard in literature) definition given in term of the spectral measure of the positive self--adjoint operator $f^*f$.
\begin{lem}\label{LL}
In the situation above $$F_{\Lambda}(f,\mu)=\operatorname{tr}_{\Lambda}\chi_{[0,\mu^2]}(f^*f)=\operatorname{dim}_{\Lambda}\operatorname{range}(\chi_{[0,\mu^2]}(f^*f))$$
as a projection in $\operatorname{End}_{\Lambda}(U).$

\noindent Notice that since $f^*f$ is a positive operator $\chi_{[0,\mu^2]}(f^*f)=\chi_{(-\ty,\mu^2]}(f^*f)$ is the spectral projection associated to the spectral resolution $f^*f=\int_{-\ty}^{\ty}\mu d\chi_{(-\ty,\mu]}.$
\end{lem}
\begin{proof}
The spectral Theorem ( a parametrized measurable version) shows that the ranges of the family of projections $\chi_{[0,\mu^2]}(f^*f)$ belong to the class $\mathcal{L}(f,\mu)$, then $$\operatorname{dim}_{\Lambda}(\operatorname{range}(\chi_{[0,\mu^2)}(f^*f)))\leq F_{\Lambda}(f,\mu).$$ In fact it's clear that $\chi_{[0,\mu^2)}(f^*_xf_x)\omega =\omega \Rightarrow \|f \omega \|\leq \mu \|\omega \|$. But for every $L\in \mathcal{L}(f,\mu)$ we get a family of injections 
$\chi_{\mu^2}(f_x^*f_x)_{|L_x}\longrightarrow \operatorname{range}(\chi_{\mu^2}(f^*_xf_x))$. Reducing modulo $\Lambda$ and using the crucial property a) of the formal dimension (Proposition 19 pag. 55 in \cite{Cos}) we get $$\operatorname{dim}_{\Lambda}(L)\leq \operatorname{dim}_{\Lambda}(\operatorname{range}(\chi_{[0,\mu^2]}(f^*f)).$$
\end{proof}

\begin{dfn}\label{elll}
A complex of random Hilbert cochains as $({L^2}(\Omega^{\bullet}X_0),d)$ and its relative and boundary versions is said {$\Lambda$--(left)--Fredholm} in degree $k$ if the differential induced on the quotient
$$\xymatrix{{\dfrac{\mathcal{D}(d^k)}{\overline{\operatorname{range}(d^{k-1})   }      }}
 \ar[r]^-d & {
 L^2(\Omega^{k+1}X_0)
 }       }$$
 gives by $\Lambda$ a.e. reduction a {left Fredholm} unbounded operator in the sense of definition \ref{frett}. In particular it happens if and only if
 \begin{equation}\label{messaggio}
 F_{\Lambda}\big{(}d|:\mathcal{D}(d^k)\cap \operatorname{range}(d^{k-1})^{\bot}\longrightarrow L^2(\Omega^{k+1}X_0),\mu\big{)}<\ty
 \end{equation} for some positive number $\mu$.

\noindent For this reason one calls the left hand--side of \eqref{messaggio} 
$$F_{\Lambda}\Big{(}L^2(\Omega^kX_0,\partial X_0),\mu\Big{)}:=F_{\Lambda}\big{(}d|:\mathcal{D}(d^k)\cap \operatorname{range}(d^{k-1})^{\bot}\longrightarrow L^2(\Omega^{k+1}X_0),\mu\big{)}$$
 the {spectral density function} of the complex at the point $k$.
\end{dfn}
\begin{rem}The definition above combined with lemma \ref{LL} says that we have to compute the formal dimension of $\chi_{[0,\mu^2]}(f^*f)$ where $f=d_{|\mathcal{D}(d)\cap \overline{\operatorname{range}(d^{k-1}})^{\bot}}.$ But as $f$ is an injective restriction of $d^k;$ every spectral projection $\chi_B(f^*f)$ projects onto a subspace that's orthogonal to $\operatorname{ker}(d^k)$. This means
\begin{equation}\label{amn}F_{\Lambda}\big{(}d|:\mathcal{D}(d^k)\cap \operatorname{range}(d^{k-1})^{\bot}\longrightarrow L^2(\Omega^{k+1}X_0),\mu\big{)}=\sup \mathcal{L}_{\Lambda}^{\bot}(f,\mu)\end{equation} where $\mathcal{L}_{\Lambda}^{\bot}(f,\mu)$ is the set of Random fields of subspaces of $\mathcal{D}(d)\cap \operatorname{ker}(d)^{\bot}$ where $d$ is bounded by $\mu$ (see Definition \ref{frett} )
\end{rem}
\noindent Now return to the boundary foliation $\mathcal{F}_{\partial}$
with its equivalence relation $\underline{\mathcal{R}_{0}}$.
\begin{theorem}
All the three complexes of Random Hilbert spaces
$$L^{2,\mathcal{F}}(\Omega^\bullet X_0),\quad
L^{2,\mathcal{F}}(\Omega^\bullet X_0,\partial X_0),\quad 
L^{2,\mathcal{F}}(\Omega^\bullet \partial X_0),$$
 with unbounded differentials are $\Lambda$--Fredholm as representations of $\underline{\mathcal{R}}=\mathcal{R}({\mathcal{F}_\partial}_0).$
\end{theorem}\label{longfre}
\begin{proof}
The proof follows by an accurate inspection of the relation between the differentials (with or without boundary conditions) and the Laplace operator trough the theory of selfadjoint boundary differential problems developed in \cite{scht}. 
To make the notation lighter let $M=L_x^0$ be the generic leaf with boundary $\partial M=\partial L_x^0$. We concentrate on the relative sequence at the point $d:A^{k}(M,\partial M)\longrightarrow A^{k+1}(M,\partial M)$ where the differential is an unbounded operator on $L^2$ with Dirichlet boundary conditions.
Let $\mathcal{D}(d)=A^{k+1}(M,\partial M)$. The following Lemma is inspired by Lemma 5.11 in \cite{lus} where in contrast Neumann boundary conditions are imposed.
\begin{lem}\label{1117}Let $\operatorname{ker}(d)$ be the kernel of $d$ as unbounded operator with Dirichlet boundary conditions, then
$$\mathcal{D}(d)\cap \operatorname{ker}(d)^{\bot}=H^1_{\textrm{Dir}}\cap \overline{\delta^{k+1} C^{\ty}_0(\Lambda^{k+1}T^*M)}^{L^2}$$ where $H^1_{\textrm{Dir}}$ is the space of order $1$ Sobolev $k$--forms $\omega$ such that $\omega_{|\partial M}=0$.
\end{lem}
\begin{proof}
First of all remember that the differential operator $d+\delta:C^{\ty}(\Lambda^{\bullet}T^*M)\longrightarrow C^{\ty}(\Lambda^{\bullet}T^*M)$ with either Dirichlet or Neumann boundary conditions is uniformly elliptic and formally self--adjoint with respect to the greenian formula
$$(d^r \omega,\eta)-(\omega,\delta^{p+1}\eta)=\int_{\partial M}(\omega \wedge \ast^{p+1}\eta)_{|\partial}.$$ 
In particular this is an elliptic boundary value problem in the classical sense according to the original definition of Lopatinski and Shapiro, together with a uniform condition on the local fundamental solutions\cite{scht}. Now let $\omega \in C^{\ty}_0$ and $\eta \in \operatorname{ker}(d)$ i.e. $\eta_n\in C^{\ty}_0$, $(\eta_n)_{|\partial M}=0$, 
 $\xymatrix{{\eta_n}\ar[r]^-{L^2}&\eta}$, 
 $\xymatrix{{d\eta_n}\ar[r]^-{L^2}&0}$ then
 $$(\eta,\delta \omega)=\lim_n(\eta_n,\delta \omega)=\underbrace{\lim_n(d\eta_n,\omega)}_{0}\pm \underbrace{\int_{\partial M}(\eta_n\wedge \ast \omega)_{|\partial M}}_{\eta_{|\partial M}=0}=0,$$ showing that 
 $\overline{\delta C^{\ty}_0}\subset \mathcal{D}(d)\cap \operatorname{ker}(d)^{\bot}.$ For the reverse inclusion take 
 
 \noindent $\omega \in \dd \cap \kerd^{\bot}$ i.e. $\omega_n\in C^{\ty}_0$, $\xymatrix{{\omega_n}\ar[r]^-{L^2}&\omega}$, $\xymatrix{{d\omega_n}\ar[r]^-{L^2}&0}.$ For fixed $\eta \in C^{\ty}_0$,
 $$\underbrace{((d+\delta)\eta,\omega)}_{d\eta \in \kerd, \omega \in \kerd^{\bot}}=(\delta \eta,\omega)=\lim_n(\delta \eta,\omega_n)\underbrace{=}_{{\omega_n}_{|\partial M}=0}=\lim_n (\eta,d\omega).$$
 \noindent Then we can apply the adjoint regularity theorem of H{\"o}rmander \cite{scht} Lemma 4.19, cor 4.22 saying that $\omega \in H^1_{\textrm{loc}}.$ It follows that $(\delta \omega,\eta)=(\omega,d \eta)$ holds because for every $\eta \in C^{\ty}_0(M-\partial M)$, $d\eta \in \kerd$ then $\delta \omega=0$. Then for every $\sigma \in C^{\ty}_0$
 $$0\underbrace{=}_{d\sigma \in \kerd}(d\sigma,\omega)=\underbrace{(\sigma,\delta \omega)}_{0}\pm \int_{\partial M}(\sigma \wedge \ast \omega)_{|\partial M}=\pm \int_{\partial M}(\overline{\omega}\wedge \overline{\ast \sigma })_{|\partial M}.$$ The last passage coming from the definition of the Hodge $\ast$ operator, 
 
 \noindent $\sigma \wedge \ast \omega=(\sigma,\omega)dvol=(\overline{\omega},\overline{\sigma})dvol=\overline{\omega}\wedge \overline{\ast \sigma}$, where $\overline{\cdot}$ is the complex conjugate in $\Lambda T^*M \otimes \mathbb {C}.$ Now from the density of $\{i^*(\overline{\ast \sigma})\}_{\sigma \in C^{\ty}_0}$ in $L^2(\partial M)$, $i:\partial M\hookrightarrow M$ the boundary condition $\omega_{|\partial M}=0$ follows. In particular $\omega \in H^1_{\textrm{Dir}}$. Now it remains to apply the Hodge decomposition  
$$L^2(\Lambda^k T^*M)=\mathcal{H}^k_{(2)}(M,\partial M)\oplus \overline{d^{k-1}\Omega_{d}^{k-1}(M,\partial M)}^{L^2}\oplus \overline{\delta^{k+1}\underbrace{\Omega_{\delta}^{k+1}(M,\partial M)}_{
\mbox{no }\partial-\mbox{conditions}
}}^{L^2} $$ to deduce $\omega \in \overline{\delta^{k+1} C^{\ty}_0(\Lambda^{k+1}T^*M)}^{L^2}.$ \end{proof}
\noindent Consider again the formally selfadjoint boundary value problem $d+\delta$ with Dirichlet boundary conditions i.e $\mathcal{D}(d+\delta)=H^1_{\mbox{Dir}}$. Its square in the sense of unbounded operators on $L^2$ is the laplacian $\Delta$ with domain 
$$H^2_{\mbox{Dir}}:=\{\omega \in H^2: \omega_{|\partial M}=0, \,((d+\delta)\omega)_{|\partial M}=(\delta \omega)_{|\partial M}=0\}.$$
\noindent Let $\Delta_k^{\bot}$ the operator obtained from $\Delta$ on $k$--forms restricted to the orthogonal complement of its kernel, it is easy to see that the splitting
$$L^2(\Lambda^k T^*M)=\mathcal{H}^k_{(2)}(M,\partial M)\oplus \overline{d^{k-1}\Omega_{d}^{k-1}(M,\partial M)}^{L^2}\oplus \overline{\delta^{k+1}\underbrace{\Omega_{\delta}^{k+1}(M,\partial M)}_{
\mbox{no }\partial-\mbox{conditions}
}}^{L^2} $$
induces the following splitting on $\Delta_k$,
$$\Delta_k^{\bot}=(\delta^{k+1}d^p)_{|\overline{\delta^{k+1}\Omega_{\delta}^{k+1}}}\oplus (d^{k-1}\delta^k)_{|\overline{d^{k-1}\Omega^{k-1}_d}}.$$
\begin{lem}\label{1118}The following identies of unbounded operators hold

 \begin{align*}&(\delta^{k+1}d^p)_{|\overline{\delta^{k+1}\Omega_{\delta}^{k+1}}}=(d^k_{|\overline{\delta^{k+1}\Omega_{\delta}^{k+1}}})^*(d^k_{|\overline{\delta^{k+1}\Omega_{\delta}^{k+1}}}),\\
&(d^{k-1}\delta^k)_{|\overline{d^{k-1}\Omega^{k-1}_d}}=(d^{k-1}_{|\overline{\delta^{k}\Omega_{\delta}^{k}}})(d^{k-1}_{|\overline{\delta^{k}\Omega_{\delta}^{k}}})^*\end{align*}
where the $d^k_{|\overline{\delta^{k+1}\Omega_{\delta}^{k+1}}}$ is the unbounded operator on the subspace $\overline{\delta^{k+1}\Omega_{\delta}^{k+1}}$ of $L^2$ with domain $H^1_{\textrm{Dir}}\cap \overline{\delta^{k+1}\Omega_{\delta}^{k+1}}$ and range $\overline{d^{k+1}\Omega^{k
+1}_d}$.
\end{lem}
\begin{proof}
This is again the \emph{dual} (in the sense of boundary conditions) statement of Lemma 5.16 in \cite{lus}.
We first state that the Hilbert space adjoint of the operator $d^k$ with domain $H^1_{\textrm{Dir}}\cap \overline{\delta^{k+1}\Omega_{\delta}^{k+1}}$ and range $\overline{d^{k+1}\Omega^{k+1}_d}$ is exactly $\delta^{k+1}$ with domain $H^1_{\textrm{Dir}}\cap \overline{d^{k}\Omega_d^{k}}$. We shall omit degrees of forms and call $d$ this restricted operator. Thanks to the intersection with $H^1$ this is also the restriction of $d+\delta$ to the same subspace, in particular $\omega \in \mathcal{D}(d^*)\subset \overline{dC^{\ty}_0}$ implies $\omega \in \dd$ and $d\omega=0.$ Take arbitrary $\eta \in H^1_{\textrm{Dir}}\cap \overline{\delta C^{\ty}_0}$, then since $\delta \eta=0$, $((d+\delta)\eta,\omega)=(d\eta,\omega)=(\eta,d^*\omega)$ and if $\eta \in H^1_{\textrm{Dir}}\cap \overline{d\Omega_d}$, $((d+\delta)\eta,\omega)=(\delta \eta,\omega)=0.$ Since $\delta H^1_{\textrm{Dir}}\bot \overline{d\Omega_d}$ this is immediately checked, 

\noindent $\sigma \in d\Omega_d$, $\sigma=d\lambda$, $\lambda_{|\partial M}=0$, $(\sigma,\delta \gamma)=\underbrace{(d\sigma,\gamma)}_{=0}+\int_{|\partial M}\underbrace{(\sigma \wedge \ast \gamma)_{|\partial M}}_{=0}.$

\noindent Also $(\eta,d^* \omega)=0$ since $d^*\omega \in \overline{\delta C^{\ty}_0}$ and $d\Omega_{\textrm{Dir}}\bot \delta C^{\ty}_0$. Then we can apply again the adjoint regularity theorem \cite{scht}, Lemma 4.19 to deduce $\omega \in H^1_{\textrm{loc}}$.
The next goal is to show $\omega \in H^1_{\textrm{Dir}}$ i.e. $d\omega,\delta \omega \in L^2$, $\omega_{|\partial M}=0$ but $dx=0\in L^2$, $\delta \omega=(d+\delta)\omega=d^* \omega \in L^2$ and

\noindent $(\omega, d\delta \eta)=(d^*\omega,\delta \eta)=(\delta \omega,\delta \eta)=(\omega,d\delta \eta)\pm \overline{\int_{\partial M}(\delta \eta \wedge \ast \omega)_{|\partial M}}$ for every $\eta \in C^{\ty}_0$. Then $0=\int_{\partial M}(\delta \eta \wedge \ast \omega)_{|\partial M}=\int_{\partial M}(\bar{\omega}\wedge \overline{\ast \delta \eta})_{|\partial M}\underbrace{=}_{=0}\int_{\partial M}({\omega}\wedge {\ast \delta \eta})_{|\partial M}$ for every $\eta$. The boundary condition follows by density.
\noindent Finally it is clear that $\delta d_{|{\mathcal{D}(d^*d)}}=\Delta=\Delta^{\bot}$ but it remains to prove the coincidence of the domains
$$\mathcal{D}(\Delta)\cap \overline{\delta C^{\ty}_0}=\mathcal{D}(d^*(d_{|\overline{\delta C^{\ty}_0}})).$$ Now $\mathcal{D}(\Delta)=H^2_{\textrm{Dir}}=\{\omega \in H^2,\,\omega_{|\partial M},\,(\delta \omega)_{|\partial M}=0\}\subset \mathcal{D}(d^*d_{|\overline{\delta C^{\ty}_0}}).$ Clearly $$\omega \in \mathcal{D}(d^*d_{|\overline{\delta C^{\ty}_0}})\Rightarrow \omega \in H^1_{\textrm{Dir}}\cap \overline{\delta C^{\ty}_0},$$ $d\omega \in H^1_{\textrm{Dir}}$ then $(d+\delta)\omega \in H^1$ and since $\omega_{|\partial M}=0$ by elliptic regularity (for the boundary value problem $(d+\delta)$ with Dirichlet conditions \cite{scht}) $\omega \in H^2$. We have just checked the boundary conditions, finally $\omega\in H^2_{\textrm{Dir}}=\mathcal{D}(\Delta)$. The second equality in the statement is proven in a very similar way.
\end{proof}
\noindent Now that the relation of $d$ with Dirichlet boundary condition restricted to the complement of its kernel with the Laplacian ($\Delta^{\bot}$) is clear we can use elliptic regularity to deduce that the relative Random Hilbert complex is $\Lambda$--Fredholm. 
This has to be done in two steps, the first is to show that the spectral function of the Laplacian controls the spectral function of the complex according to the equation
\begin{equation}\label{fdg}
F_{\Lambda}(\Delta^{\bot}_k,\sqrt{\mu})=F_{\Lambda}(L^{2,\mathcal{F}}(\Omega^k X_0,\partial X_0),\mu)+F_{\Lambda}(L^{2,\mathcal{F}}(\Omega^{k-1} X_0,\partial X_0),\mu).
\end{equation}
In fact \begin{align}\nonumber F_{\Lambda}(\Delta^{\bot}_k,\sqrt{\mu})&=F_{\Lambda}\Big{(}(\delta^{k+1}d^k)_{|\overline{\delta^{k+1}\Omega^{k+1}_{\delta}}  }),\sqrt{\mu}\Big{)}+F_{\Lambda}\Big{(}(d^{k-1}{\delta}^k)_{|\overline{d^{k-1}\Omega^{k-1}_{d}}  }),\sqrt{\mu}\Big{)}\\ \nonumber
&=
F_{\Lambda}\Big{(}(d^k_{|\overline{\delta^{k+1}\Omega_{\delta}^{k+1}}})^*(d^k_{|\overline{\delta^{k+1}\Omega_{\delta}^{k+1}}}),\sqrt{\mu}\Big{)}+F_{\Lambda}\Big{(}(d^{k-1}_{|\overline{\delta^{k}\Omega_{\delta}^{k}}})(d^{k-1}_{|\overline{\delta^{k}\Omega_{\delta}^{k}}})^*,\sqrt{\mu}\Big{)}\\ \nonumber
&=F_{\Lambda}\Big{(}d^{k}_{|\overline{\delta^{k+1}\Omega_{\delta}^{k+1}}},\mu\Big{)}+F_{\Lambda}(d^{k-1}_{|\delta^k\Omega_{\delta}^k},\mu\Big{)}.
\end{align}In the first step we have used the obvious fact that the spectral functions behave additively under direct sum of operators togheter with the remark after \eqref{elll}. At the second step there are lemmas \ref{1117} and \ref{1118} together with the following properties of the spectral functions. 
\begin{itemize}
\item $F_{\Lambda}\Big{(}f^*f,\sqrt{\lambda}\Big{)}=F_{\Lambda}(f,\lambda)$
\item $F_{\Lambda}(\phi,\lambda)=F_{\lambda}(\phi^*,\lambda)$
\end{itemize} These can be adapted 
to hold in our situation with unbounded operators. Good references are the paper of Lott and L{\"u}ck \cite{lolu} and the paper of L{\"u}ck and Schick \cite{lus}.

\noindent Let us firs recall the equation
\begin{equation}\nonumber
F_{\Lambda}(\Delta^{\bot}_k,\sqrt{\mu})=F_{\Lambda}(L^{2,\mathcal{F}}(\Omega^k X_0,\partial X_0),\mu)+F_{\Lambda}(L^{2,\mathcal{F}}(\Omega^{k-1} X_0,\partial X_0),\mu). 
\end{equation}
It says that
 we have only to show that $\Delta_k^{\bot}$ is left $\Lambda$--Fredholm to have control of both Fredholmness at degree $k$ and $k-1$. We can use the heat kernel, in fact by elliptic regularity for each leaf the heat kernel 
$e^{-t {\Delta_k,x}^{\bot}}(z,z')$ is smooth and uniformly bounded along the leaf on intervals $[t_0,\ty)$ \cite{scht} Theorem 2.35. As $x$ varies in $\partial X_0$ these bounds can made uniform by the uniform geometry (in fact the constants depend on the metric tensor, its inverse and a finite number of their derivatives in normal coordinates) and we get a family of smooth kernels that varies transversally in a measurable fashion since it is obtained by functional calculus from a measurable family of operators. Then they give a $\Lambda$--trace class element in the relevant Von neumann algebra. Now the projections $\chi_{[0,\mu]}(f^*f)$ in definition \ref{elll} where $f$ is the differential restricted to the complement of its kernel are obtained from the heat kernel as $$\chi_{[0,\mu]}(f^*f)=\underbrace{\chi_{[0,\mu]}(\Delta_k^{\bot})e^{\Delta_k^{\bot}}}_{\textrm{bounded}}
\underbrace{\chi_{[0,\mu]}(\Delta_k^{\bot})e^{-\Delta_k^{\bot}}}_{\Lambda-\textrm{trace class}}<\ty.$$
\end{proof}
\begin{rem}
The same argument of elliptic regularity for b.v. problems together with the various Hodge decompositions shows that each term of the long sequence \eqref{success} {is a finite} {Random Hilbert space.}
\end{rem}
\begin{theorem}The long sequence \eqref{success}
\begin{equation*}
\xymatrix{{\cdot \cdot \cdot}\ar[r]& H^{k,\mathcal{F}}_{dR,(2)}(X_0,\partial X_0)\ar[r]^-{i_{*}} & 
H^{k,\mathcal{F}}_{dR,(2)}(X_0)\ar[r]^-{r_{*}}&\\
{}&\ar[r]^-{r_{*}}&
H^{k,\mathcal{F}}_{dR,(2)}(\partial X_0)\ar[r]^-{\delta}&
H^{k-1,\mathcal{F}}_{dR,(2)}(X_0,\partial X_0)\ar[r]& \cdot \cdot \cdot}
\end{equation*} is $\Lambda$--weakly exact (definition \ref{vonexact} )
\end{theorem}
\begin{proof}
This is exactly the same proof of Cheeger and Gromov \cite{ChGr} ( see also the book by L{\"u}ck 
\cite{luk}
(Theorem 1.21). In fact the crucial final step there, that is based on the property of formal dimension of Hilbert $\Gamma$--modules
$$\operatorname{dim}_{\Gamma}\Big{(}\bigcap_{i\in I}V_i\Big{)}=\inf_{i\in I} \operatorname{dim}_{\Gamma} V_i,$$ can be replaced here by the corresponding property for the formal dimension of Random Hilbert spaces (a proof in the appendix)
\end{proof}

\subsubsection{The proof}
\begin{theorem}We have
$$\sigma_{\Lambda,dR}(X_0,\partial X_0)=\sigma_{\Lambda,\textrm{an}}(X,\partial X_0)$$ thus together with formula \eqref{anasign}  all the three signatures we have defined agree
$$\sigma_{\Lambda,dR}(X_0,\partial X_0)=\sigma_{\Lambda,\operatorname{an}}(X_0,\partial X_0)=\sigma_{\Lambda,\operatorname{Hodge}}(X_0,\partial X_0)=\langle L(X),C_{\Lambda}\rangle +1/2[\eta_{\Lambda}(D^{\mathcal{F}_{\partial}})]$$
\end{theorem}
\begin{proof}
We pass through different intermediate results, sometimes doing leafwise considerations. Our model is of course the work of L{\"u}ck and Schick \cite{lus}. The proof of L{\"u}ck and Schick in turn is inspired by the classical argument of Atiyah Patodi and Singer \cite {AtPaSi1} with the issue that at $L^2$ level long sequences are only weakly exact and the spectrum of the boundary operator is not discrete.
\bigskip

\noindent {\bf{First step}}. This is done. We have proved the equality $$\sigma_{\Lambda,\operatorname{an}}(X_0,\partial X_0)=\sigma_{\Lambda,\operatorname{Hodge}}(X_0,\partial X_0)$$ where at right--hand side there is the signature on harmonic leafwise $L^2$--forms on the elonged manifold with elongated foliation i.e. the $\Lambda$ signature of the Poincar\`e product on leafwise harmonic forms. Our reference in this section is then the harmonic signature. 
\bigskip

\noindent {\bf{Second step}}.
We shall prove $\sigma_{\Lambda,dR}(X_0,\partial X_0)=\sigma_{\Lambda,\operatorname{Hodge}}(X_0,\partial X_0)$. We explain now the strategy

We have to measure the $+/-$ eigenspaces of the intersection form on the field of Hilbert spaces $H^{2k}_{dR,(2)}(X_0,\partial X_0)$ as square integrable representations of 
$\mathcal{R}_0$ (the whole foliation on $X_0$). Now thanks to the fundamental note on section \ref{fundationtr} it is sufficient to measure the corresponding projections in the von Neumann algebra arising by restriction of the Random Hilbert spaces to 
$\underline{\mathcal{R}_0}$ (the equivalence relation of the foliation induced on the boundary). This is a consequence of the very definition of the trace as an integral of a functor with values measure spaces and the fact the boundary contains a complete transversal. The passage to 
$\underline{\mathcal{R}_0}$ has the great vantage we can write boundary problems and sequences of random Hilbert spaces, in particular the third term  in 
$$\xymatrix{
0\ar[r]&A_x^{k-1}(L_x^0,\partial L_x^0)\ar[r]^-i&A_x^{k-1}(L_x^0)\ar[r]^-r &A_x^{k-1}(\partial L_x^0)\ar [r]&0
}$$ as representations of $\underline{\mathcal{R}_0}.$

\noindent Remember the notation: $x\in \partial X_0$, $L_x^0$ is the leaf of the compact foliated manifold with boundary, $L_x$ is the leaf of the foliation on the manifold $X$ with a cylinder attached.
Consider the random Hilbert space $H^{2k}_{dr,(2)}(X_0)$ obtained from the various $L^2$ cohomologies of the leaves with no boundary conditions (this is called sometimes in literature the $L^2$--homology since it naturally pairs with forms with Dirichlet boundary conditions). We have a family of restrictions $\partial X_0\ni x \longmapsto r_x^p:\mathcal{H}^{2k}(L_x)\longrightarrow H^{2k}_{dR,(2)}(L_x^0)$ and intertwining operators   $(\mathcal{H}^{2k}(L_x))_{x\in X_0}:\longmapsto H^{2k}_{dR,(2)}(L_x^0)$. There are also natural mappings $i_x^{2k}:H^{2k}_{dR,(2)}(L_x^0,\dl)\longrightarrow H^{2k}_{dR,(2)}(L_x)$ and the mappings $q$ coming from the long sequence in cohomology \begin{equation}\label{prty}\xymatrix{H^{2k}_{dR,(2)}(L_x^0,\dl)\ar[r]^-{i_x^{2k}}&H^{2k}_{dR,(2)}(L_x^0)\ar[d]^-{q^{2k}_x}\\
\mathcal{H}^{2k}(L_x)\ar[ur]^-{r_x^{2k}}&               H^{2k}_{dR,(2)}(\dl)}.
\end{equation}

\noindent Following the program of L{\"u}ck and Schick we shall prove 
\begin{enumerate}

\item  
{\begin{itemize}
\item \begin{equation}\label{gess}
\overline{\operatorname{range}(r^{2k})}=\overline{\operatorname{range}(i^{2k})}\textrm{ as projections in } \operatorname{End}_{\Lambda,\underline{\mathcal{R}_0}}\Big{[}H_{dR,(2)}^{2k,\mathcal{F}}(X_0)\Big{]}.\end{equation}  
\item  The signature can be computed looking at the fields of sesquilinear Poincar\`e products on the images of $i_x^{2k}$ as square integrable representations of $\mathcal{R}_0$,
\begin{equation}\label{closure}\xymatrix{H^{2k}_{dR,(2)}(L_x^0,\dl)\ar[r]^-{i_x^{2k}}&H^{2k}_{dR,(2)}(L_x^0)\\
\mathcal{H}^{2k}(L_x)\ar[ur]^-{r_x^{2k}}&{}                           }\quad \quad  x\in \mathcal{R}_0.\end{equation}
\end{itemize}}

\item The signature of the field of products on the image of $i_x^{2k}$ concides with the signature of the fields of Poincar\'e products on $(\mathcal{H}_x)_{x\in X_0}$ as square integrable representations of $\underline{\mathcal{R}_0}$ {\bf{that in turn coincides with those computed tracing in}} $\mathcal{R}_0$
\end{enumerate}
\noindent Notice about \eqref{gess} that $\overline{\operatorname{range}(i^{2k})}=\ker q^{2k}$ by the long exact sequence.
\bigskip

\noindent 1.
\noindent L{\"u}ck and Schick (lemma 3.12 in \cite{lus}) prove the following result 
\begin{lem}
$$q_x^{2k} \circ r_x^{2k}=0,\,\, x\in \partial X_0$$
\end{lem}\label{lafine}
 \noindent By definition of the algebra of intertwining operators $$q_x^{2k}\circ r^{2k}_x=0\,\forall x\in \partial X_0 \Longrightarrow \overline{\ra r_x^{2k}}\subset \ker q_x^2k$$ then  $\overline{\ker q^{2k}}\cdot \overline{\ra r^{2k}}=\overline{\ra i^{2k}}\cdot \overline{\ra r^2k}=\overline{\ra r^{2k}}\in \operatorname{End}_{\Lambda,\underline{\mathcal{R}_0}}\Big{[}H_{dR,(2)}(X_0)\Big{]}.$
 
\noindent Now Von Neumann dimentions come in play in a fundamental way. Consider the field of unbounded boundary differentials $d_x:L^2(\Omega^{2k-1}\partial L_x^0)\longrightarrow L^2(\Omega^{2k-1}\partial L_x^0)$ exactly as in \cite{lus} (and essentially by elliptic regularity and the fact trace=trace on the boundary foliation) they define a left Fredholm affiliated operator so the image of the field of the spectral projection $\chi_{(0,\gamma]}(\delta d)$ has dimension going to zero for $\gamma\rightarrow 0$. Given $\epsilon>0$ define the following field of subspaces,
$$E_{\epsilon,x}^{2k}:=\ra (d\circ \chi_{(\gamma,\ty)}(\delta d ))\subset L^2(\Omega^{2k}\partial L_x^0).$$ Properties of $E_\epsilon^{2k}:$
\begin{enumerate}
\item \underline{it is measurable}, in fact is obtained by functional calculus from a natural Borel family.
\item \underline{It has codimension less that} $\epsilon$ in $\overline{\ra(d)}$ in fact $d\circ \chi_{(-\ty,0]}(\delta d)=0.$
\item \underline{It is closed}, because the restriction of $\delta d$ to the subspace corresponding to $(0,\ty)$ satisfies $\delta d \geq \gamma$ than is invertible (seen using the polar decomposition).
\end{enumerate}
\noindent Now we have to invoke the leafwise Hodge decomposition with (Neumann)  boundary condition, 
\begin{equation}
\label{hokka}
L^2(\Omega^{2k-1}(L_x^0))=\overline{\ra d^{2k-2}}\oplus \overline{ \ra \delta^{2k-2}_{|\{\omega_{|\partial}=0\}}}\oplus 
 \ker \Delta^{2k}_{|\{ (\ast \omega)_{|\partial}=0=(\delta \omega)_{|\partial}      \}}.
\end{equation} 
The methods of Schick \cite{scht}  surely applies to the generic leaf $L_x^0$ in fact this is bounded geometry and has a collar so the fact its boundary has infinite connected components (complete in the induced metric) plays no role. So the space $\clis$ can be canonically identified with the third addendum in \eqref{hokka} and pull back to the boundary gives a well defined measurable
\footnote{inverse image of a measurable field of subspaces by a unif. bounded measurable family of bounded operators is measurable, one can split the domain space as $\operatorname{Ker}\oplus\operatorname{Ker}^{\bot}$ and apply the well known fact that inverses of isom. are measurable \cite{Dix}} 
family of (uniformely) bounded mappings $\beta_x^{2k}:\clis \longrightarrow \cali$. Define, by pull--back the following measurable field of closed subspaces
$$K_{\epsilon,x}^{2k}\subset \clis.$$
\noindent Properties of $K_{\epsilon,x}^{2k}:$
\begin{enumerate}
\item $K_{\epsilon,x}^{2k}\subset \clis$
\item $K_{\epsilon,x}^{2k}\subset (\beta^{2k}_x)^{-1}(\overline{\ra d_{\partial}})$
\item The field $K_{\epsilon,x}^{2k}$ defines a projection having codimension in $\ker q^{2k}$ that's less than $\epsilon.$
\end{enumerate}
\noindent Then there's another density lemma in \cite{lus} (Lemma 3.16) 
stating 
$$K^{2k}_{\epsilon}\subset \ra(r^{2k}_x:\mathcal{H}^{2k}_{2}(L_x)\longrightarrow \overline{\ra i^{2k}}).$$
\noindent All of this properties certainly say that \eqref{gess} is true (by normality of the trace we can reach $\overline{\ra (i^{2k})}$ with a family of subprojections whose codimension tends to zero).
\bigskip

\noindent 2. 

\noindent Again following \cite{lus}, $q[0]^{2k}$ (notation of the proof of Lemma \ref{lafine}) defines a bounded family of mappings from $\mathcal{H}^{2k}_{(2)}(L_x)$ to $\overline{\ra d_{\partial}}$. So let $\mathcal{H}_{\epsilon,x}^{2k}\subset \mathcal{H}^{2k}_{(2)}(L_x)$ be as before the inverse image of $E_{\epsilon,x}^{2k}$. Since we are using harmonic forms the pull--back is (uniformly) bounded in the $L^2$ norm so $\mathcal{H}_{\epsilon,x}^{2k}$ is a field of closed subspaces giving projection of codimension in $\mathcal{H}^{2k}_{(2)}(X)$ not greater than $\epsilon$. Now if $L_{\epsilon,x}^{2k}\subset \overline{\ra i^{2k}_x}$ is the closure of the image of $L_{\epsilon,x}^{2k}$ under the mapping $r^{2k}_x:\mathcal{H}^{2k}_{(2)}(L_x)\longrightarrow \overline{\ra i^{2k}_x},$ its codimension into $\overline{\ra i^{2k}}$ is less than $\epsilon$ exactly because of \eqref{gess} since the codimension of $\mathcal{H}^{2k}_\epsilon$ in $\mathcal{H}^{2k}_{(2)}(X)$ is less than $\epsilon$.

\noindent The leafwise intersection form $$s_x^0:H^{2k}_{(2)}(L_x^0,\partial L_x^0)\times H^{2k}_{(2)}(L_x^0,\partial L_x^0)\longrightarrow \C$$ descends into a pairing on $\overline{\ra i^{2k}_x}$  which restricts to 
$\eta_x^0:L_{\epsilon,x}^{2k}\times L_{\epsilon,x}^{2k}\longrightarrow \C.$
\noindent The codimension of $L_{\epsilon}^{2k}\subset \overline{\ra i^{2k}}$ is less than $\epsilon$ one gets
$$|\operatorname{sign}_{\Lambda}(s^0)-\operatorname{sign}_{\Lambda}(\eta)|\leq \epsilon,$$
remember that $\operatorname{sign}_{\Lambda}(s^0)=\sigma_{\Lambda,dR}(X_0,\partial X_0).$

\noindent Now it's a quite amazing computation performed by Luck and Schick \cite{lus} that the leafwise Hodge intersection form we called $s^{\ty}_x:\mathcal{H}_{(2)}(L_x)\times \mathcal{H}_{(2)}(L_x)\longrightarrow \C$ \underline{descends} to a pairing on each $\mathcal{H}_{\epsilon,x}^{2k}$ and in turn to exactly the pairing $\eta_x^0$ defined above. Again since the codimension of $\mathcal{H}_{\epsilon}^{2k}$ in $\mathcal{H}^{2k}_{(2)}(X)$ is $\leq \epsilon$ we get $|\operatorname{sign}_{\Lambda}(s^{\ty})-\operatorname{sign}_{\Lambda}(\eta^0)|\leq \epsilon$ then
$$|\operatorname{sign}_{\Lambda}(s^\ty)-\operatorname{sign}_{\Lambda}(s^0)|\leq 2\epsilon.$$ The theorem is proved since $\epsilon$ is arbitrary. \end{proof}

\begin{rem}{\bf{On the assumption of the complete transversal contained into the boundary.}}
The assumption $\operatorname{Saturation}(\partial X_0)=X_0$ is really simply 
avoidable in fact one can write the sequence
$$\xymatrix{
0\ar[r]&A_x^{k-1}(L_x^0,\partial L_x^0)\ar[r]^-i&A_x^{k-1}(L_x^0)\ar[r]^-r &A_x^{k-1}(\partial L_x^0)\ar [r]&0
}$$
for $x$ also in the interior but the last arrow is null for $\partial L_x^0=0$ so everything works in the exactly same way. Otherwise divede $X_0$ into the Borel sets made of leaves touching/not touching the boundary and proceed exactly as above.
\end{rem}

\begin{rem}{\bf{Foliated $\rho$--invariants.}}
The $L^2$ signature formula of Vaillant \cite{Vai} has been used by Piazza and Schick in \cite{PS} to prove the existence  of a large class of non trivial $L^2$--rho invariants. Now $\rho$ invariants for foliated flat bundles are defined by Benameur and Piazza in \cite{BP} but the definition works on general (measured) foliations. The author believe that a similar (more sophisticated) tecnique can be reproduced to show that formula \eqref{ffo} implies the existence of non trivial foliated rho--invariants. It should be stressed that up to now no examples are known.
\end{rem}
\appendix
\section{Von Neumann algebras, foliations and index theory}\label{appe}

\subsection{Non--commutative integration theory.}\label{paoletto}
The measure--theoretical framework of non--commutative integration theory is particular fruitful when applied to measured foliations. 
\noindent The non--commutative integration theory of Alain Connes \cite{Cos} provides us a measure theory on every measurable groupoid $(G,\mathcal{B})$ with $G^{(0)}$ the space of unities. In our applications $G$ will be mostly the equivalence relation $\mathcal{R}$ or sometimes the holonomy groupoid of a foliation. Transverse measures in the non--commutative integration theory sense will be defined from holonomy invariant transverse measures. Below a list of fundamental objects and facts. This is a very brief and simplified survey in fact the general theory admits the existence of a modular function that says, in the case of foliations how the transverse measure of sets changes under holonomy (under flows generated by fields tangent to the foliation). Hereafter our modular function is everywhere 1, corresponding to the geometrical case of a foliation equipped with a \underline{holonomy invariant transverse measure} (this is a definition we give below).

\begin{description}
\item[Measurable groupoids]. A groupoid is a small category $G$ where every arrow is invertible. The set of objects is denoted by $G^{(0)}$ and there are two maps 
$s,r:G\longrightarrow G^{(0)}$ where $\gamma:s(\gamma)\longrightarrow r(\gamma).$ Two arrows 
$\gamma_1,\gamma_2$ can be composed if 
$r(\gamma_2)=s(\gamma_1)$ and the result is $\gamma_1\cdot \gamma_2$. The set of composable arrows is $G^{(2)}=\{(\gamma_1,\gamma_2):r(\gamma_2)=s(\gamma_1)\}$. As a notation 
$G_x=r^{-1}(x)$,  $G^x=s^{-1}(x)$ for $x\in G^{(0)}.$
 An equivalence relation $\mathcal{R}\subset X\times X$ is a groupoid with $r(x,y)=x$ and $s(x,y)=y$, in this manner $$(z,x)\cdot (x,y)=(z,y).$$ The range of the map $(r,s):G\longrightarrow G^{(0)}\times G^{(0)} $ is an equivalence relation called {\underline{the principal groupoid associated}} to $G$. In this sense groupoids desingularize equivalence relations.
A \underline{measurable groupoid} is a pair 
$(G,\mathcal{B})$ where $G$ is a groupoid and 
$\mathcal{B}$ is a $\sigma$--field on $G$ making measurable the structure maps $r,s$, composition 
$\circ:G^{(2)}\longrightarrow G$ and the inversion $\gamma \longmapsto \gamma^{-1}$. 
\item[Kernels] are mappings $x\longmapsto \lambda^x$ where $\lambda^x$ is a positive measure on $G$, supported on the $r$--fiber $G^x=r^{-1}(x)$ with a measurability property i.e. for every set $A\in \mathcal{B}$ the function $y\longmapsto \lambda^{y}(A)\in [0,+\infty]$ must be measurable.

\noindent
A kernel $\lambda$ is called {\underline{proper}} if there exists an increasing family of measurable sets $(A_n)_{n\in \mathbb{N}}$ with $G=\cup_n A_n$ making the functions $\gamma\longmapsto \lambda^{s(\gamma)}(\gamma^{-1}(A))$ bounded for every $n\in \mathbb{N}$. The point here is that every element $\gamma:x\longrightarrow y$ in $G$ defines by left traslation a measure space isomorphism $G^x\longrightarrow G^y$ and calling \begin{equation}\label{horml} R(\lambda)_{\gamma}:=\gamma \lambda^x\end{equation} (here $\gamma \lambda_x$ is push--forward measure under the $\gamma$--right traslation) one has a kernel in the usual sense i.e. a mapping with value measures. The definition of properness is in fact properness for $R(\lambda)$. 

\noindent The space of proper kernels is denoted by $\mathcal{C}^+$.
\item[Transverse functions] are kernels 
$(\nu^{x})_{x\in X}$ with the left invariance property 
$$\gamma \nu^{s(\gamma)}=\nu^{r(\gamma)}$$ for every 
$\gamma \in G.$ One checks at once that properness is equivalent to the existence of an increasing family of measurable sets 
$(A_n)_n$ with $G=\cup_n A_n$ such that the functions 
$x\longmapsto \nu^{x}(A_n)$ are bounded for every $n\in \mathbb{N}$. The space of proper transverse functions is denoted $\mathcal{E}^{+}.$ 

\noindent The {\underline{support}} of a transverse function $\nu$ is the measurable set $$\operatorname{supp}(\nu)=\{x\in G^{(0)}:\nu^x\neq 0\}.$$ This is $\underline{saturated}$ w.r.t. the equivalence relation induced by $G$ on $G^{(0)}$, $ x\mathcal{R} y$ iff there exists $\gamma:x\longrightarrow y$. If $\operatorname{supp}(\nu) =G^{(0)}$ we say that $\nu$ is {\underline{faithful}}. 

\noindent When $G=\mathcal{R}$ or the holonomy groupoid this gives families of positive measures one for each leaf in fact in the first case the invariance property is trivial, in the second case we are giving a measure $\nu^x$ on each holonomy cover $G^x$ with base point $x$ but the invariance property says that these are invariant under the deck trasformations together with the change of base points then push forward on the leaf under $r:G^x\longrightarrow L_x.$ 
\item[Convolution.] The groupoid structure provides an operation on kernels. For fixed kernels $\lambda_1$ and, $\lambda_2$ on $G$ their convolution product is the kernel $\lambda_1\ast\lambda_2$ defined by
$$(\lambda_1 \ast\lambda_2)^y=\int(\gamma \lambda_2^x)d\lambda_1^y(\gamma),\quad y\in X.$$ It is a fact that if $\lambda$ is a kernel and $\nu$ is a transverse function then $\nu \ast \lambda$ is a transverse function. Clearly $R(\lambda_1\ast \lambda_2)=R(\lambda_1)\circ R(\lambda)$ the standard composition of kernels on a measure space. Here $R(\cdot)$ is that of equation \eqref{horml}.
\item[Transverse invariant measures] \label{tras}(actually transverse measures of modulo $\delta=1$). These are linear mappings $\Lambda:\mathcal{E}^+\longrightarrow [0,+\infty]$ such that
{\begin{enumerate}
\item $\Lambda$ is normal i.e 
$\Lambda(\sup \nu_n)=\sup \Lambda(\nu_n)$ for every increasing sequence 
$\nu_n$ in $\mathcal{E}^+$  bounded by a transverse function. Since the sequence is bounded by an element of $\mathcal{E}^+$ the expression $\sup \nu_n$ makes sense in $\mathcal{E}^+$.
\item $\Lambda$ is invariant under the right traslation of $G$ on $\mathcal{E}^+$. This means that $$\Lambda(\nu)=\Lambda(\nu \ast \lambda)$$ for every $\nu\in \mathcal{E}^+$ and kernel $\lambda$ such that $\lambda^y(1)=1$ for every $y\in G^{(0)}$.
\end{enumerate}}
\noindent A transverse measure is called {\underline{semi--finite}} if it is determined by its finite values i.e $$\Lambda(\nu)=\sup\{\Lambda(\nu'),\, \nu'\leq \nu,\, \Lambda(\nu')<\infty\}.$$ We shall consider only semi--finite measures.

\noindent A transverse measure is 
$\sigma$--{\underline{finite}} if there exists a faithful transverse function 
$\nu$ of kind $\nu=\sup{\nu_n}$ with $\lambda(\nu_n)<\infty$.

\noindent The coupling of a transverse function $\nu\in \mathcal{E}^+$ and a transverse measure $\Lambda$ produces a positive measure $\Lambda_{\nu}$ on $G^{(0)}$ through the equation $\Lambda_{\nu}(f):=\Lambda((f\circ s)\nu.$ The invariance property reflects downstairs (in the base of the groupoid) in the property $\Lambda_{\nu}(\lambda)=\Lambda(\nu \ast \lambda)$ for $\nu \in \mathcal{E}^+$ and $\lambda \in \mathcal{C}^+$.

\noindent Measures on the base $G^{(0)}$ that can be represented as $\Lambda_{\nu}$ are characterized by a theorem of disintegration of measures.

\begin{theorem}\label{pocahontas}(Connes \cite{Cos})
Let $\nu$ be a transverse proper function with support $A$. 

\noindent The mapping $\Lambda \longmapsto \Lambda_{\nu}$ is a bijection between the set of transverse measures on the reduced groupoid $$G_A^A=r^{-1}(A) \cup s^{-1}(A)$$ and the set of positive measures $\mu$ on $G^{(0)}$ satisfying the following equivalent relations
\begin{enumerate}
\item $(\mu \circ \nu)^{\tilde{}}
=\mu \circ \nu$
\item $\lambda,\lambda'\in \mathcal{C}^+, \nu \ast \lambda=\nu \ast \lambda' \in \mathcal{\epsilon}^+ \Longrightarrow \mu(\lambda(1))=\mu(\lambda'(1)).$
\end{enumerate}
\end{theorem}
Nex we shall explain this procedure of disintegration in a geometrical way for foliations.

\noindent We shall see that what is important here is the 
class of null--measure subsets of $G^{(0)}$. A saturated set $A\subset G^{(0)}$ is called $\Lambda$--{\underline{negligible}} if $\Lambda_{\nu}(A)=0$ for every $\nu\in \mathcal{E}^{+}$.

\item[Representations.] Let $H$ be a measurable field of Hilbert spaces; by definition this is a mapping $x\longmapsto H_x$ from $G^{(0)}$ with values Hilbert spaces. The measurability structure is assigned by a linear subspace $\mathcal{M}$ of the free product vector space of the whole family $ \Pi_{x\in G^{(0)}}H_x$ meaning that
\begin{enumerate}
\item For every $\xi \in \mathcal{M}$ the function $x\longmapsto \|\xi(x)\|$ is measurable.
\item A section $\eta \in  \Pi_{x\in G^{(0)}}H_x$ belongs to $\mathcal{M}$ if and only if the function $\langle \eta(x),\xi(x)\rangle$ is measurable for every $\xi \in \mathcal{M}$.
\item There exists a sequence $\{\xi_i\}_{i\in \mathbb{N}}\subset \mathcal{M}$ such that $\{\xi_i(x)\}_{i\in \mathbb{N}}\subset \mathcal{M}$ is dense in $H_x$ for every $x$. 
\end{enumerate} Elements of $\mathcal{M}$ are called measurable sections of $H$.

\noindent Suppose a measure $\mu$ on $G^{(0)}$ has been chosen. One can put together the Hilbert spaces $H_x$ taking their direct integral
$$\int H_x d\mu(x).$$ This is defined as follows, first select the set of square integrable sections in $\mathcal{M}$. This is the set of sections $s$ such that the integral $\int_{G^{(0)}}\|s(x)\|_{H_x}^2d\mu(x)<\ty$ then identify two square integrable sections if they are equal outside a $\mu$--null set. The direct integral comes equipped with a natural Hilbert space structure with product $$\langle s,t\rangle:=\int_{G^{(0)}} \langle s(x),t(x)\rangle_{H_x}d\mu(x).$$ The notation $s=\int_{G^{(0)}}s(x)d\mu(x)$ for an element of the direct integral is clear. A field of bounded operators $x\longmapsto B_x \in B(H_x)$ is called \underline{measurable} if sends measurable sections to measurable sections. A mesurable family of operators with operator norms uniformely bounded i.e. $\operatorname{ess} \sup \|B_x\|<\ty$ defines a bounded operator called \underline{decomposable} $B:=\int_{G^{(0)}}B_xd\mu(x)$ on the direct integral in the simplest way $$Bs:=\int_{G^{(0)}}B_xd\mu(x)\, s=\int_{G^{(0)}}B_xs(x)d\mu(x).$$ For example each element of the abelian Von Neumann algebra $L^{\ty}_{\mu}(G^{(0)})$ defines a decomposable operator acting by pointwise multiplication.
One gets an involutive algebraic isomorphism of $L^{\ty}_{\mu}(G^{(0)})$ onto its image in $B(\int H_x d\mu(x))$ called the algebra of \underline{
diagonal operators}.
 One can ask when a bounded operator $T\in B(\int H_x d\mu(x))$ is decomposable i.e. when $T=\int T_x d\mu(x)$ for a family of uniformely bounded operators $(T_x)_x$.
The answer is precisely when it belongs to the commutant of the diagonal algebra.

\noindent A {\underline{representation}} of $G$ on $H$ is the datum of an Hilbert space isomorphism $U(\gamma):H_{s(\gamma)}\longrightarrow H_{r(\gamma)}$ for every $\gamma\in G$ with 
\begin{enumerate}
\item $U(\gamma_1^{-1}\gamma_2)=U(\gamma_1)^{-1}U(\gamma_2),\quad \forall \gamma_1,\gamma_2\in G,\quad r(\gamma_1)=r(\gamma_2).$
\item For every couple $\xi,\eta$ of measurable section the function defined on $G$ according to
$\gamma \longmapsto \langle \eta_{r(\gamma)},U(\gamma)\eta_{s(\gamma)}\rangle,$ is measurable.
\end{enumerate}
\noindent A fundamental example is given by the {\underline{left regular representation}} of $G$ defined by a proper transverse function $\nu \in \mathcal{E}^+$ in the following way. The measurable field of Hilbert space is $L^2(G,\nu)$ defined by $x\longmapsto L^2(G^x,\nu^x)$ with the unique measurable structure making measurable the family of sections of the kind $y\longmapsto f_{|G^x}$ obtained from every measurable $f$ on $G$ such that each $\int |f|^2d\nu^x$ is finite. For every $\gamma:x\longrightarrow y$ in $G$ one has the left traslation $L(\gamma):L^2(G^x,\nu^x)\longrightarrow L^2(G^y,\nu^y)$,\quad $(L(\gamma)f)(\gamma')=f(\gamma^{-1}\gamma')$, $\gamma'\in G^y$. 

\item[Intertwining operators] are morphisms between representations. If $(H,U)$, $(H',U')$ are representations of $G$ an intertwining operator is a measurable family of operators $(T_x)_{x\in G^{(0)}}$ of bounded operators $T_x:H_x\longrightarrow H'_x$ such that
\begin{enumerate}
\item Uniform boundedness; $\operatorname{ess}\sup \|T_x\|<\infty$.
\item For every $\gamma\in G$  $U'(\gamma)T_{s(\gamma)}=T_{r(\gamma)}U(\gamma).$
\end{enumerate}
\noindent Looking at a representation as a measurable functor, an intertwining operator gives a natural transformation between representations. The vector space of intertwing operators from $H$ to $H'$ is denoted by $\operatorname{Hom}_{G}(H,H').$
\item[Square integrable representations.] Fix some transverse function $\nu \in \mathcal{E}^+$. For a representation of $G$ the property of being equivalent to some sub--representation of the infinite sum of the regular left representation $L^{\nu}$ is independent of $\nu$ and is the definition of {\underline{square integrability}} for representations. This is a generalization of the concept of square integrable representations of locally compact groups. Actually, due to measurability issues much care is needed here to define sub representations (see section 4 in \cite{Cos}) but the next fundamental remark assures that square integrable representations are very commons in applications.  

\item[Measurable functors and  representations.] Let $\tilde{\mathcal{R}}_+$ be the category of (standard) measure spaces without atoms i.e. objects are triples $(\mathcal{Z},\mathcal{A},\alpha)$ where $(\mathcal{Z},\mathcal{A})$ is a standard measure space and $\alpha$ is a $\sigma$--finite positive measure while morphisms are measurable mappings.  

\noindent Measurability of a functor $F:G\longrightarrow \tilde{\mathcal{R}}_+$ is a measure structure on the disjoint union $Y=\bigcup_{x\in G^{(0)}}F(x)$ making the following structural mappings measurable
\begin{enumerate}
\item The projection $\pi:Y\longrightarrow G^{(0)}.$
\item The natural bijection $\pi^{-1}(x)\longrightarrow F(x).$
\item The map $x\longmapsto \alpha^x$, a $\sigma$--finite measure on $F(x)$.
\item The map sending $(\gamma,z)\in G\times X$ with $s(\gamma)=\pi(z)$ into $F(\gamma)z\in Y$.
\end{enumerate} 
\noindent Usually one assumes that $Y$ is union of a denumerable collection $(Y_n)_n$ making every function $\alpha^x(Y_n)$ bounded. With a measurable functor $F$ one has an associated representation of $G$ denoted by $L^2\bullet F$ defined in the following way: the field of Hilbert space is $x\longmapsto L^2(F(x),\alpha^x)$ and if $\gamma:x\longrightarrow y$ then define $U(\gamma):L^2(F(x),\alpha^x)\longrightarrow L^2(F(y),\alpha^y)$ by $f\longmapsto F(\gamma^{-1})\circ f.$ Proposition 20 in \cite{Cos} shows that this is a square--integrable representation.
\item[Random hilbert spaces and Von Neumann algebras.] We have seen that every fixed transverse measure $\Lambda$ defines a notion of $\Lambda$--null measure sets (for saturated sets) hence an equivalence relation on $\operatorname{End}_{G}(H_1,H_2)$
the vector space of all intertwining operators $T:H_1\longrightarrow H_2$ between two square integrable representations $H_i$. Each equivalence class is called a {\underline{random operator}} and the set of random operators is denoted by $\operatorname{End}_{\Lambda}(H_1,H_2)$. 
\noindent Also square integrable representations can be identified modulo $\Lambda$--null sets. An equivalence class of square integrable representations is by definition a {\underline{random hilbert space}}.  
\noindent 

\noindent Theorem 2 in \cite{Cos} says that $\operatorname{End}_{\Lambda}(H)$ is a Von Neumann algebra for every random Hilbert space. 

\noindent More precisely  choose some $\nu \in \mathcal{E}^+$ and put $\mu=\Lambda_{\nu}$ and $m:=\mu \bullet \nu$ to form the Hilbert space $\mathcal{H}=L^2(G,m)$. For a function $f$ on $G$ denote $Jf=f^{\sharp}(\gamma)=\bar{f(\gamma^{-1})}$, consider the space $\mathcal{A}$ of measurable functions $f$ on $G$ such that $f,f^{\sharp}\in L^2(G,m)$ and $\sup(\nu|f^{\sharp}|)<\infty$. Equip $\mathcal{A}$ with the product $f\ast_{\nu}g=f\nu\ast g$. The structure $\mathcal{A}$ has is that of an {\underline{Hilbert algebra}} (a left--Hilbert algebra in the modular case) i.e $\mathcal{A}$ is a $\ast$--algebra with positive definite (separeble) pre--Hilbert structure such that
\begin{enumerate}
\item $\langle x,y\rangle=\langle y^*,x^*\rangle, \quad \forall x,y \in \mathcal{A}.$
\item The representation of $\mathcal{A}$ on $\mathcal{A}$ by left multiplication is bounded, involutive and faithful.
\end{enumerate}
\noindent With such structure one can speak about the left regular representation $\lambda$ of 
$\mathcal{A}$ on the Hilbert space completion $\mathcal{H}$ of $\mathcal{A}$ itself. The double commutant 
$\lambda''(\mathcal{A})$ of this representation is the Von Neumann algebra $W(\mathcal{A})$ associated to the Hilbert algebra 
$\mathcal{A}$. It is a remarkable fact that 
$W(\mathcal{A})$ comes equipped with a semifinite faithful normal trace $\tau$ such that 
$$\tau(\lambda(y^*)\tau(x))=\langle x,y\rangle \quad \forall x,y \in \mathcal{A}.$$
\noindent Furthermore one knows that the commutant of $\lambda(\mathcal{A})$ in ${\mathcal{H}}$ is generated by the algebra of right multiplications $\lambda'(\mathcal{A})=J\lambda(\mathcal{A})J$ for the conjugate--linear isometry $J:\mathcal{H}\longrightarrow \mathcal{H}$ defined by the involution in $\mathcal{A}$.
\noindent For every $\Lambda$--random Hilbert space $H$ one can use the measure $\Lambda_{\nu}$ on $G^{(0)}$ to form the {\underline{direct integral}} $\nu(H)=\int H_xd\Lambda_{\nu}(x)$. Remember that the direct integral is the set of equivalence classes modulo $\Lambda_{\nu}$ zero measure of square integrable measurable sections.
 Now, directly from the definition, an intertwining operator $T\in \operatorname{Hom}_{\Lambda}(H_1,H_2)$ is a decomposable operator defining a bounded operator $\nu(T):\nu(H_1)\longrightarrow \nu(H_2).$

Put $W(\nu)$ for the Von Neumann algebra associated to the Hilbert algebra $L^2(G,m)$, $m=\Lambda_{\nu}\bullet \nu$, $\nu\in \mathcal{E}^+$. 
\begin{theorem}\label{equivalenza}(Connes)
Fix some transverse function $\nu \in \mathcal{E}^+$
\begin{enumerate}
\item For every $\Lambda$--random Hilbert space $H$ there exists a unique normal representation of $W(\nu)$ in $\nu(H)$ such that $U_{\nu}(f)=U(f\nu)$ $f\in \mathcal{A}_{\nu}.$ Here $U(f\nu)$ is defined by $(U(f\nu)\xi)_y=\int U(\gamma)\xi_xd(f\nu^y)(\gamma).$
\item The correspondence $H\longmapsto \nu(H)$, $T\longmapsto \nu(T)$ is a functor from the $(W^*)$--category $\mathcal{C}_{\Lambda}$ of random Hilbert spaces and intertwining operators to the category of $W(\nu)$ modules.
\item If the transverse measure $\nu$ is faithful the functor above is an equivalence of categories.
\end{enumerate}
\end{theorem}
\noindent Then in the case of faithful transverse measures one gets an isometry of $\operatorname{End}_{\Lambda}(H)$ on the commutant of $W(\nu)$ on the direct integral $\nu(H)$. In particular $\operatorname{End}_{\Lambda}(H)$ is a Von Neumann algebra.
\item[Transverse integrals.] The most important notion of non commutative integration theory is the integral of a {\underline{random variable}} against a transverse measure. A positive random variable on $(G,\mathcal{B},\Lambda)$ is nothing but a measurable functor $F$ as defined above. Let $$X:=\bigcup_{x\in G^{(0)}}F(x)$$ be the disjoint union measure space and $\bar{\mathcal{F}}^+$ the space of measurable functions with values in $[0,+\infty]$ while ${\mathcal{F}}^+$ 
is for functions with values on 
$(0,+\ty]$. Kernels $\lambda$ on $G$ acts as convolution kernels on  
$\bar{\mathcal{F}}^+$ according to $$(\lambda \ast f)(z)=\int f(\gamma^{-1}z)d\lambda^{y}(\gamma), \quad 
y=\pi(z)\in G^{(0)}.$$ This is an associative operation $(\lambda_1\ast \lambda_2)\ast f=\lambda_1 \ast (\lambda_2\ast f)$.

\noindent Now to define the integral $\int Fd\lambda$ choose some faithful $\nu$  and put
$$\int F d\lambda=\sup\{\Lambda_{\nu}(\alpha(f)),\, f\in {\mathcal{F}}^+,\,\nu \ast f\leq 1\},$$ 
this is independent from $\nu$ and enjoys the following properties \begin{enumerate}
\item there exist random variables $F_1,F_2$ with $F=F_1\oplus F_2$ such that $\int F_1d\Lambda=0$ and a function $f_2\in \mathcal{F}^+(X_2)$ with $X_2=\bigcup_{x\in G^{(0)}}F_2(x)$ with $\nu \ast f_2=1.$
\item Monotony. If $f,f'\in \mathcal{F}(X)$ satisfy $\nu \ast f\leq \nu \ast f'\leq 1$ then $$\Lambda_{\nu}((\alpha(f))\leq \Lambda_{\nu}((\alpha(f'))$$ in particular for $F_2$ as in 1.
$$\int F_2 d\Lambda=\Lambda_{\nu}((\alpha(f')).$$
\end{enumerate}
\item[Traces.] 
Let $A$ be a Von Neumann algebra with the cone of positive elements $A^+$. 

\noindent A \underline{weight} on a $A$ is a functional $\phi:A^+\longrightarrow [0,\infty]$ such that
\begin{enumerate}
\item $\phi(a+b)=\phi(a)+\phi(b)$, $a,b\in A^+$
\item $\phi(\alpha a)=\alpha \phi(a)$, $\alpha \in \R^+$, $a\in A^+$.
\end{enumerate}
a weight is a 
\underline{trace} if it is invariant under interior automorphisms, put in another way
$$\phi(a^*a)=\phi(aa^*),\,a \in A^+.$$ 
A weight is called 
\begin{itemize}
\item \underline{faithful} if $\phi(a)=0\Rightarrow a=0$, $a\in A^+$.
\item \underline{normal} if for every increasing net $\{a_i\}_i$ of positive elements with least upper bound $a$ then
$$\phi(a)=\sup \{\phi(a_i)\}.$$
\item \underline{Semifinite} if the linear span of a the set of $\phi$--finite elements, $\{a\in A^+:\phi(a)<\infty\}$ is $\sigma$--weak dense. 
\end{itemize} Every V.N algebra has a semifinite normal faithful weight. This is not true for traces, V.N algebras that can be equipped with a semifinite faithful trace are called \underline{semifinite}. 
The Von Neumann algebra $\operatorname{End}_{\Lambda}(H)$ associated to a square integrable representation of a measurable groupoid is always semifinite (this is an effect of square integrability) and comes equipped with a bijection $T\longmapsto \Phi_T$ between positive operators and semifinite normal weights $\Phi_T:\operatorname{End}_{\Lambda}(H)\longrightarrow [0,+\ty]$ where $\Phi_T$ is faithful if and only if $T_x$ is not singular $\Lambda$--a.e. The construction of this correspondence uses the fact, that for a faithful transverse function $\nu$ the
direct integral $\nu(H)=\int H_x d \Lambda_{\nu}(x)$ is a module over the
 Von Neumann algebra $W(\nu)$ associated to the Hilbert algebra $\mathcal{A}$ above described. 
 
 \noindent The notation of Connes is $$\Phi_T(1):=\int \operatorname{Trace}(T_x)d\Lambda(x)$$ i.e. the mapping $T\longmapsto \Phi_T(1)$ is the canonical trace on $\operatorname{End}_{\Lambda}(H)$. In fact this is related to the type $I$ Von Neumann algebra 
 $P$ of classes modulo equality $\Lambda_{\nu}$ almost everywhere of measurable fields $(B_x)_{x\in G^{(0)}},\,B_x\in B(H_x)$ of bounded operators. Remember that $P$ has a canonical trace $\rho(B)=\int \operatorname{Trace}(B_x)d\Lambda_{\nu}(x)$ hence we can define
 $$\rho_T(B):=\int \operatorname{Trace}(T_xB_x)d\Lambda_{\nu}(x).$$
 \begin{lem}\label{11010}
For a faithful transverse function $\nu$ there's a unique operator valued weight\footnote{see \cite{take} for the definition} $E_{\nu}$ from $P$ to $\operatorname{End}_{\Lambda}(H)$ such that the diagram
$$
\xymatrix{{P^+}\ar[d]_{E_{\nu}}\ar[dr]^{\,\,\,\,\,\,\,\,\,\,\,\rho_T(\cdot)=\int\operatorname{Trace }(T_x\cdot)d\Lambda_{\nu}(x)} & \\
\operatorname{End}_{\Lambda}(H)\ar[r]_{\Phi_T} & \mathbb{C}}
$$is commutative. Moreover $E_{\nu}$ is such that
$$E_{\nu}(B)=C$$
 if $B=(B_x)_{x\in G^{(0)}}$, $B\in P^+$ is an operator making bounded the corresponding family $$C_y:=\int U(\gamma)B_xU(\gamma)^{-1}d\nu^y.$$
\end{lem}

Let $F$ be a random variable and put $H=L^2\bullet F$. The integration process above defines a semi--finite faithful trace on the Von Neumann algebra $\operatorname{End}_{\Lambda}(H)$. In fact, for $T\in\operatorname{End}^{+}_{\Lambda}(H)$ let ${F_T}$ the new random variable defined by $x\mapsto(F(x),\alpha_T(x))$ where $\alpha_T(x)$ is the measure on $F(x)$ such that $$\alpha_T(x)(f)=\operatorname{Trace_{L^2}}(T_x^{1/2}M(f)T_x^{1/2})$$ where $f$ is a bounded measurable function on $F(x)$ and $M(f)$ the corresponding multiplication operator on $L^2(F(x))$. The trace is 
$$\Phi_T(1)=\int F_Td \Lambda.$$ 

\noindent In the following we shall use often the notation $\operatorname{tr}_{\Lambda}(T)=\Phi_T(1)$ to emphasize the dependence on $\Lambda$.

\noindent With a trace one can develop a dimesion theory for square integrable representation i.e. a dimension theory for random Hilbert spaces that's very similar to the dimension theory of $\Gamma$--Hilbert modules.

\bigskip

\noindent The \underline{formal dimension} of the random Hilbert space $H$ is $$\operatorname{dim}_{\Lambda}(H)=\int \operatorname{Trace}(1_{H_x})d\Lambda(x)$$ here some fundamental properties
\begin{lem}\label{formaldimension}
\begin{enumerate}
\item  If $\operatorname{Hom}_{\Lambda}(H_1,H_2)$ contains an invertible element then $\operatorname{dim}_{\Lambda}(H_1)=\operatorname{dim}_{\Lambda}(H_2).$
\item $\operatorname{dim}_{\Lambda}(\oplus H_i)=\sum \operatorname{dim}_{\Lambda}(H_i).$
\item $\operatorname{dim}_{\sum \Lambda_i}(\oplus H)=\sum \operatorname{dim}_{\Lambda_i}(H).$
\end{enumerate}
 \end{lem}
 \item[Formal dimensions and projections]
  We need more properties of the formal dimension that are implicit in Connes work but not listed above.
 
\noindent Start to consider sub--square integrable representation.
Consider a Random Hilbert space $(H,U)$; if for every $x$ one choose in a mesurable way a closed subspace $K$ such that $U(\gamma):K_x\longrightarrow K_y$ for every $\gamma \in G$ we say that $(K,V)$, $V(\gamma):=U(\gamma)_{|i_x K_x}$ is a \underline{sub Random Hilbert space.} Once a faithful $\nu\in \mathcal{E}^+$ is keeped fixed, the functor $\nu$ in theorem \ref{equivalenza} displays $H$ and $K$ as submodules of the V.N. algebra associated to the Hilbert Algebra $\mathcal{A}$, hence there must be an injection $\operatorname{End}_{\Lambda}(K)\longrightarrow \operatorname{End}_{\Lambda}(H)$. In fact from the diagram 
$$
\xymatrix{{W(\nu)}\ar[r] \ar[dr]& B(\nu(H)) & \\
{}&B(\nu(K))\ar[u]_{\nu(i)}
}
$$
we see that multiplication by the bounded operator $\nu(i)=\int_{G^{(0)}}\nu(i_x)d\Lambda_{\nu}(x)$ sends the commutator of $W(\nu)$ in $B(\nu(K))$ into its commutator in $B(\nu(H))$. To check that the natural traces $\varphi^H\in P(\operatorname{End}_{\Lambda}(H))$ and $\varphi^K\in P(\operatorname{End}_{\Lambda}(K))$ are preserved by this inclusion we can examine with more detail the meaning of square integrability for a representation. So let us consider the subset of measurable vector fields
$$D(V,\nu):=\Big{\{} \xi:\exists c>0: \forall y \in G^{(0)},\forall \alpha \in K_y
\int_{G^y}|\langle \alpha, V(\gamma)\xi_x\rangle_{K_y} |^2d\nu(\gamma)\leq c^2 \|\alpha\|^2
 \Big{\}}.$$ The definition of square integrability is equivalent to the statement that $D(V,\nu)$  contains a denumerable total subset. 
 In other words the operation of assigning a coefficient $\alpha\longmapsto T_{\nu}(\xi)\alpha=(\alpha,\xi)$ defines an intertwining operator from $V$ to the left regular representation of $L^{\nu}$ of $G$, on the field of Hilbert spaces $L^2(R^x,\nu^x)_x$. This has the property $T_{\nu}(\xi)^*f=V(f\nu)\xi,\quad \xi\in D(V,\nu)$ if $\nu|f|$ is bounded. Then, for $\xi,\eta \in D(V,\nu)$ the operator 
 $$\theta_{\nu}(\xi,\eta):=T_{\nu}(\xi)^*T_{\nu}(\eta) \in \operatorname{End}_{\Lambda}(K)$$ 
 satisfies the following interesting identity $(\theta_{\nu}(\xi,\eta)\xi',\eta')=(\xi',\eta)\ast_{\nu}(\eta',\xi)^{v}$ for bounded measurable sections of $K$. Furthermore the vector space $\mathcal{J}_{\nu}$ generated by couples $\xi,\eta \in D$ is a \underline{bilateral ideal} and respects ordering for transverse functions, 
 $$\mathcal{J}_{\nu}\subset \mathcal{J}_{\nu'} \textrm{ if } \nu\leq \nu'.$$ Since the measure $\nu$ is faithful this is also \underline{weakly dense} hence completely determines the trace by the simple formula
 \begin{equation}\label{novembre}
 \varphi^K(\theta_{\nu}(\xi,\xi))=\int_{G^{(0)}}\langle \xi_x,\xi_x\rangle d\Lambda_{\nu}(x),\quad \xi \in D(V,K).
 \end{equation} 
Now via $i$ 
we get an inclusion $D(V,\nu)\subset D(U,\nu)$ let's check this statement: let $\xi\in D(V,\nu)$, $y\in G^{(0)}$, $\alpha \in H_y$ then
\begin{align*}
\int_{G^y}|\langle \alpha,U(\gamma)i_x\xi_x\rangle_{K_x}|^2d\nu(\gamma)&=\int_{G^y}|\langle \alpha,i_x V(\gamma)\xi_x\rangle |_{K_x}^2d\nu(\gamma)\\
&=
\int_{G^y}|\langle P_{H_y}\alpha,V(\gamma)\xi_x\rangle|^2d\nu(\gamma)\leq c^2 \|\alpha\|^2.
\end{align*} 
It turns out that under the inclusion $\operatorname{End}_{\Lambda}(K)\hookrightarrow \operatorname{End}_{\Lambda}(H)$ it is essential to check how a $\theta_{\nu}(\xi,\xi)$ acts on $H$ and to check that the two natural traces are equal. These two problems are very simple. Since for $\xi \in D(V,K)$ the endomorphism $\theta_{\nu}(\xi,\xi)$ under the inclusion is sent in $\operatorname{End}_{\Lambda}(H)$ to the operator
 \begin{equation}\label{proppp}
 \theta_{\nu}(i\xi,i\xi)=T_{\nu}(i\xi)^*T(i\xi)=i^*T_{\nu}(\xi)T_{\nu}(\xi)i
 \end{equation}
  we can prove the following
 \begin{lem}\label{comppppa}
\begin{enumerate}
\item The natural traces are compatible w.r.t. the inclusions, in other words we have a commutative diagram
$$
\xymatrix{{\operatorname{End}_{\Lambda}(K)}\ar[d]^{\varphi^K}\ar[r]&\operatorname{End}_{\Lambda}(H)\ar[dl]^{\varphi^H}\\
\R
}
$$
\item To get the formal dimension of $K$ as a Random Hilbert space is sufficient to trace the corresponding field of projections in $\operatorname{End}_{\Lambda}(H)$
\end{enumerate}
 \end{lem}
 \begin{proof}By the computation \eqref{proppp} above, the density result on the ideal $\mathcal{J}_{\nu}$ and formula \eqref{novembre}
 it is suffcient to check the next identity
 \begin{align*}\varphi^H(\theta_{\nu}(i\xi,i\xi))&=
 \varphi^H(T_{\nu}(i\xi)^*T_{\nu}(i\xi))=
 \int_{G^{(0)}}\langle i_x\xi_x,i_x\xi_x\rangle_{H}d\Lambda_{\nu}(x)\\&=\int_{G^{(0)}}\langle \xi_x,\xi_x\rangle_K d\Lambda_{\nu}(x)=\varphi^K(\theta_{\nu}(\xi,\xi)).
 \end{align*}
 \end{proof}
 \noindent Now we have the tools to prove two crucial properties of the formal dimension similar to the properties of the dimension of $\Gamma$-- Hilbert modules (compare Chapter 1. of \cite{luk} )
 \begin{prop}\label{decrescenza}
 Let $\{(H^{(i)},U^{(i)})\}_{i\in I}$ a system of Random Hilbert subspaces of $(H,U)$ directed by $\subset$ then
$$\operatorname{dim}_{\Lambda}
(
\operatorname{closure}\Big{(}\bigcup_{i\in I}H^{(i)}\Big{)}
\Big{)}=
\sup\{\operatorname{dim}_{\Lambda}H_i,i\in I\}$$
\noindent If the system is directed by $\supset$ then
$$\operatorname{dim}_{\Lambda}
\Big{(}\bigcap_{i\in I}H^{(i)}
\Big{)}=
\inf\{\operatorname{dim}_{\Lambda}H_i,i\in I\}$$
\end{prop}
\begin{proof}
The choice of a faithful normal transverse function $\nu \in \mathcal{E}^+$ estabilishes the equivalence of categories described above between $C_{\Lambda}$ and the category of normal representations of the Von Neumann algebra associated with $W(\nu)$; the first statement then follows from the compatibility of the natural traces proved in \ref{comppppa} and the \underline{normality} (the passage to $\sup$) of the trace in the limit square integrable representation. 
The second statement follows from the first adopting a standard trick changing a decreasing system into an increasing one. It is in fact sufficient to consider $H^{(i)\bot}$ and observe $${\Big{(}\bigcup_{i\in I}H^{(i)\bot}\Big{)}^{\bot}}=\bigcap_{i\in I}H^{(i)}.$$ From the fact that the family is bounded by $H$ we can write the following equation with finite $\Lambda$--dimensions
$$\operatorname{dim}_{\Lambda}\Big{(}     H^{(i)\bot}       \Big{)}=\operatorname{dim}_{\Lambda}(H)-\operatorname{dim}_{\Lambda}(H^{(i)})$$
\end{proof}
 \end{description}

\end{document}